\input ./style/arxiv-general.cfg
\documentclass[MSNbibl,number,citesort,seceqn,dvips]{arxbj}
\makeatletter
   \@ifpackageloaded{graphicx}{}{\usepackage{graphicx}}
\makeatother

\volume{23}
\issue{1}
\pubyear{2017}
\firstpage{249}
\lastpage{287}
\doi{10.3150/15-BEJ743}
\docsubty{FLA}

\makeatletter
\newcommand{\rrvert}{\vert}
\newcommand{\llvert}{\vert}
\def\overset{\stackrel}
\newcommand{\eqref}[1]{(\ref{#1})}
\def\IND{\mathbb{I}}
\def\PROB{\mathbb{P}}
\def\EXP{\mathbb{E}}
\def\mathe{\mathrm{e}}
\def\var{\operatorname{Var}}
\def\cov{\operatorname{Cov}}
\def\N{\mathbb{N}}
\def\R{\mathbb{R}}

\newtheorem{thmm}{Theorem}[section]
\newremark{rem}{Remark}[section]

\newtheorem{cor}[thmm]{Corollary}

\newtheorem{prop}[thmm]{Proposition}

\newproclaim{dfn}{Definition}
\newremark{ex}{Example}
\makeatother

\begin{document}
\begin{frontmatter}

\title{Concentration inequalities in the infinite urn scheme for
occupancy counts and the missing mass, with applications}
\runtitle{Concentration inequalities in the infinite urn scheme}

\begin{aug}
\author[A]{\inits{A.}\fnms{Anna} \snm{Ben-Hamou}\corref{}\thanksref{A}\ead[label=e1]{anna.benhamou@univ-paris-diderot.fr}},
\author[B]{\inits{S.}\fnms{St\'ephane} \snm{Boucheron}\thanksref{A,B}\ead[label=e2]{stephane.boucheron@univ-paris-diderot.fr}\ead[label=u1,url]{stephane-v-boucheron.fr}}
\and\\
\author[C]{\inits{M.I.}\fnms{Mesrob I.} \snm{Ohannessian}\thanksref{C}\ead[label=e3]{mesrob@gmail.com}\ead
[label=u2,url]{sites.google.com/site/mesrob}}
\address[A]{LPMA UMR 7599 CNRS \& Universit\'e Paris-Diderot,
5, rue Thomas Mann 75013, Paris, France.
\printead{e1}}
\address[B]{DMA ENS Ulm,
45, rue d'Ulm 75005, Paris, France.\\
\printead{e2,u1}}
%
%
\address[C]{University of California, San Diego, Atkinson Hall, 9500
Gilman Dr, La Jolla, CA 92093, USA.
\printead{e3,u2}}
\end{aug}

%
\received{\smonth{12} \syear{2014}}
%
\revised{\smonth{5} \syear{2015}}

%
\begin{abstract}
An infinite urn scheme is defined by a probability mass function
$(p_j)_{j\geq1}$ over positive integers. A random allocation consists
of a sample of $N$ independent drawings according to this probability
distribution where $N$ may be deterministic or Poisson-distributed.
This paper is concerned with occupancy counts, that is with the number
of symbols with $r$ or at least $r$ occurrences in the sample, and with
the missing mass that is the total probability of all symbols that do
not occur in the sample. Without any further assumption on the sampling
distribution, these random quantities are shown to satisfy
Bernstein-type concentration inequalities. The variance factors in
these concentration inequalities are shown to be tight if the sampling
distribution satisfies a regular variation property. This regular
variation property reads as follows. Let the number of symbols with
probability larger than $x$ be $\vec{\nu}(x) = |\{ j \colon p_j \geq
x\}
|$. In a regularly varying urn scheme, $\vec{\nu}$ satisfies $\lim_{\tau
\rightarrow0} \vec{\nu}(\tau x)/\vec\nu(\tau) = x^{-\alpha}$ for
$\alpha\in[0,1]$ and the variance of the number of distinct symbols
in a sample tends to infinity as the sample size tends to infinity.
Among other applications, these concentration inequalities allow us to
derive tight confidence intervals for the Good--Turing estimator of the
missing mass.
\end{abstract}

%
\end{frontmatter}
%
\section{Introduction}
\label{sec:introduction}

From the 20th century to the 21st, various disciplines have tried to
infer something about scarcely observed events: zoologists about
species, cryptologists about cyphers, linguists about vocabularies, and
data scientists about almost everything. These problems are all about
`small data' within possibly much bigger data. Can we make such inference?

\textit{Problem setting}.
To move into a concrete setting, let $U_1, U_2,\ldots, U_n$ be i.i.d.
observations from a fixed but unknown distribution $(p_j)_{j=1}^{\infty
}$ over a discrete set of symbols $\mathbb{N}_{+}=\mathbb{N}\setminus
\{
0\}$. We consider each $j$ in $\mathbb{N}_{+}$ as a discrete \emph
{symbol} devoid of numerical significance. The terminology of `infinite
urn scheme' comes from the analogy to $n$ independent throws of balls
over an infinity of urns, $p_j$ being the probability of a ball falling
into urn $j$, at any $i$th throw. We alternatively adhere to the
symbols or the urns perspective, based on which carries the intuition
best. Species, cyphers, and vocabularies all being discrete, are well
modeled as such. The sample size $n$ may be fixed in advance; we call
this the \emph{binomial setting}. It may be randomly set by the
duration of an experiment; this gives rise to the \emph{Poisson
setting}. More precisely, in the latter case we write it as $N$, a
Poisson random variable independent of $(U_i)$ and with expectation
$t$. We index all Poisson-setting quantities by $t$ and write them with
functional notations, instead of subscripts used for the fixed-$n$ scheme.

For each $j,n\in\mathbb{N}_{+}$, let $X_{n,j}= \sum_{i=1}^n \mathbb
{I}_{\{U_i=j\}}$ be the number of times symbol $j$ occurs in a sample
of size $n$, and $X_j(t)=\sum_{i=1}^{N(t)} \mathbb{I}_{\{U_i=j\}}$ the
Poisson version. In questions of underrepresented data, the central
objects are sets of symbols that are repeated a small number $r$ of
times. The central quantities are the \emph{occupancy counts} $K_{n,r}$
[respectively, $K_r(t)$ for the Poisson setting], defined as the number
of symbols that appear exactly $r$ times in a sample of size $n$:
\[
K_{n,r}= \bigl\vert\{j, X_{n,j}=r\} \bigr\vert= \sum
_{j=1}^\infty\IND_{\{
X_{n,j}=r\}}.
\]

The collection $(K_{n,r})_{r\geq1}$ [resp. $(K_{r}(t))_{r\geq1}$] has
been given many names, such as the ``profile'' (in information theory
(Orlitsky \textit{et al.} \cite{MR2095850})) or the ``fingerprint'' (in
theoretical computer
science (Batu \textit{et al.} \cite{Batu2001}, Valiant and Valiant
\cite
{Valiant2011})) of the probability distribution
$(p_j)_{j \in\mathbb{N}_{+}}$. Here we refer to them by \emph
{occupancy counts} individually, and \emph{occupancy process} all together.

The occupancy counts then combine to yield the cumulated occupancy
counts $K_{n,\overline{r}}$ [respectively $K_{\overline{r}}(t)$] and
the total number of distinct symbols in the sample, or the total number
of occupied urns, often called the \emph{coverage} and denoted by $K_n$
[respectively $K(t)$]:
\[
K_{n,\overline{r}}= \bigl\vert\{j, X_{n,j}\geq r\} \bigr\vert= \sum
_{j=1}^\infty\IND_{\{X_{n,j}\geq r\}} = \sum
_{s \geq r} K_{n,s}
\]
and
\[
K_n= \bigl\vert\{j, X_{n,j}>0\} \bigr\vert= \sum
_{j=1}^\infty\IND_{\{
X_{n,j}>0\}} = \sum
_{r \geq1} K_{n,r}.
\]

In addition to the occupancy numbers and the number of distinct
symbols, we also address the \emph{rare} (or small-count) \emph
{probabilities} $M_{n,r}$ [respectively, $M_r(t)$], defined as the
probability mass corresponding to all symbols that appear exactly $r$ times:
\[
M_{n,r}= \PROB\bigl(\{j, X_{n,j}=r\}\bigr) = \sum
_{j=1}^\infty p_j \IND_{\{
X_{n,j}=r\}}.
\]
In particular, we focus on $M_{n,0}=\sum_{j=1}^\infty p_j \IND_{\{
X_{n,j}=0\}}$, which is called the missing mass, and which corresponds
to the probability of all the unseen symbols.

Explicit formulas for the moments of the occupancy counts and masses
can be derived in the binomial and Poisson settings. The occupancy
counts' expectations are given by
\begin{eqnarray*}
\EXP K_n&=&\sum_{j=1}^\infty
\bigl(1-(1-p_j)^n\bigr),\qquad \EXP K(t) =\sum
_{j=1}^\infty\bigl(1-\mathe^{-tp_j}\bigr),
\\
\EXP K_{n,r}&=& \sum_{j=1}^\infty
\pmatrix{n
\cr
r} p_j^r(1-p_j)^{n-r},\qquad
\EXP K_r(t) =\sum_{j=1}^\infty
\mathe^{-tp_j} \frac{(tp_j)^r}{r!},
\\
\EXP M_{n,r}& =& \sum_{j=1}^\infty
\pmatrix{n
\cr
r} p_j^{r+1}(1-p_j)^{n-r},\qquad
\EXP M_r(t) = \sum_{j=1}^\infty
\mathe^{-tp_j}p_j \frac{(tp_j)^r}{r!} .
\end{eqnarray*}
%
Formulas for higher moments can also be computed explicitly but their
expression, especially in the binomial setting where a lot of
dependencies are involved, often has an impractical form.

This classical occupancy setting naturally models a host of different
application areas, including computational linguistics, ecology, and
biology. Urns may represent species, and we are interested in the
number of distinct species observed in a sample (the ecological
diversity) or in the probability of the unobserved species. In
linguistics, urns may represent words. In both of these applications,
the independence assumption of the random variables $\{U_i\}
_{i=1,\ldots,
n}$ may seem unrealistic. For instance in a sentence, the probability
of appearance of a word strongly depends on the previous words, both
for grammatical and semantic reasons. Likewise, the nucleotides in a
DNA sequence do not form an i.i.d. sample. In $n$-gram models,
independence is only conditional and the observations are assumed to
satisfy a Markovian hypothesis: the probability of occurrence of a word
depends on the $n-1$ previous words. But the i.i.d. case, although very
simple, yields results that are interesting in themselves, and upon
which a more sophisticated framework may be built.

Many practical questions may now be formulated in this setting. If we
double the duration of a first experiment, how many yet unobserved
specimens would we find (how does $K_{2n,r}$ compare to $K_{n,r}$
[resp. $K_r(2t)$ to $K_r(t)$]) (Fisher \textit{et al.} \cite{Fisher1943})?
If certain cypher
keys have been observed, what is the probability for the next to be
different (how does one estimate $M_{n,0}$)? For instance, Good~\cite
{Good1953} and Turing observed that $(n+1)\mathbb{E}M_{n,0} = \mathbb
{E} K_{n+1,1}$ for all $n\geq1$, and proposed to estimate the missing
mass using the Jackknife estimator $G_{n,0} = K_{n,1}/n$ (the
proportion of symbols seen just once). 

\textit{Contributions}.
To study the Good--Turing estimator or other quantities that depend
significantly on the small-count portion of the observations, we need
to understand the occupancy counts well. Our contribution here is to
give sharp concentration inequalities with explicit and reasonable
constants, for $K_n$, $K_{n,r}$, and $M_{n,0}$ [resp. $K(t)$, $K_r(t)$,
$M_0(t)$]. We give distribution-free results, and then exhibit a vast
domain where these results are tight, namely the domain of
distributions with a heavier tail than the geometric. In this domain,
the non-asymptotic exponential concentration properties that we
establish are sharp in the sense that the exponents are order-optimal,
precisely capturing the scale of the variance. For this reason, we
dedicate a portion of the paper to establishing bounds on various variances.

\textit{Organization}.
The paper is organized as follows. In Section~\ref{sec:notation}, we
present our terminology and give a concise summary of the results. In
Section~\ref{sec:concentration} we present our variance bounds and
concentration results for the occupancy counts and the missing mass in
great generality. In Section~\ref{sec:regular-variation}, we specialize
these results to regularly varying distributions, the aforementioned
domain of distributions where concentration can be characterized
tightly. We then elaborate on some applications in Section~\ref
{sec:applications}, and conclude with a discussion of the results and
possible extensions in Section~\ref{sec:discussion}. We group all
proofs in the end, in Section~\ref{sec:proofs}.

\section{Summary of results}
\label{sec:notation}

\textit{Terminology}.
Our concentration results mostly take the form of bounds on the
log-Laplace transform. Our terminology follows closely (Boucheron
\textit{et al.}
\cite
{BoLuMa13}). We say that the random variable $Z$ is \emph{sub-Gaussian}
on the right tail (resp. on the left tail) with variance factor $v$ if,
for all $\lambda\geq0$ (resp. $\lambda\leq0$),
%
\begin{equation}
\label{eq:4} \log\EXP\mathe^{\lambda(Z-\EXP Z)} \leq\frac{v\lambda^2}{2}.
\end{equation}
%

We say that a random variable $Z$ is \emph{sub-Poisson} with variance
factor $v$ if, for all $\lambda\in\R$,
%
\begin{equation}
\label{eq:5} \log\EXP\mathe^{\lambda(Z-\EXP Z)} \leq v\phi(\lambda) ,
\end{equation}
with $\phi:\lambda\mapsto\mathe^\lambda-\lambda-1$.

We say that a random variable $Z$ is \emph{sub-gamma on the right tail}
with variance factor $v$ and scale parameter $c$ if
%
\begin{equation}
\log\EXP\mathe^{\lambda(X-\EXP X)} \leq\frac{\lambda
^2v}{2(1-c\lambda
) } \qquad\mbox{for every }\lambda
\mbox{ such that } 0\leq\lambda\leq 1/c.\label{eq:6}
\end{equation}
The random variable $Z$ is \emph{sub-gamma on the left tail} with
variance factor $v$ and scale parameter~$c$, if $-Z$ is sub-gamma on
the right tail with variance factor $v$ and scale parameter~$c$.
If $Z$ is sub-Poisson with variance factor $v$, then it is sub-Gaussian
on the left tail with variance factor $v$, and sub-gamma on the right
tail with variance factor $v$ and scale parameter $1/3$.

These log-Laplace upper bounds then imply exponential tail bounds. For
instance, inequality \eqref{eq:6} results in a Bernstein-type
inequality for the right tail, that is, for $s>0$ our inequalities have
the form
\[
\PROB\bigl\{Z > \EXP[Z]+\sqrt{2vs} + cs\bigr\} \leq\mathe^{-s},
\]
while inequality \eqref{eq:4} for all $\lambda\leq0$ entails
\[
\PROB\bigl\{Z < \EXP[Z]-\sqrt{2vs} \bigr\} \leq\mathe^{-s}.
\]
We present such results first without making distributional
assumptions, beyond the structure of those quantities themselves. These
concentrations then specialize in various settings, such as that of
regular variation.

\textit{Main results}.
We proceed by giving a coarse description of our main results. In the
Poisson setting, for each $r\geq1$, $(\mathbb{I}\{X_j(t)=r\})_{j\geq
1}$ are independent, hence $K_{r}(t)$ is a sum of independent Bernoulli
random variables, and it is not too surprising that it satisfies
sub-Poisson, also known as Bennett, inequalities. For $\lambda\in
\mathbb{R}$, we have:
\[
\log\mathbb{E} \mathe^{\lambda(K_r(t) -\mathbb{E}K_r(t))} \leq \operatorname{var}
\bigl(K_r(t)\bigr) \phi(\lambda) \leq\mathbb{E}\bigl[K_r(t)
\bigr] \phi (\lambda).
\]
The proofs are elementary and are based on the careful application of
Efron--Stein--Steele inequalities and the entropy method (Boucheron
\textit{et al.}
\cite{BoLuMa13}).

As for the binomial setting, the summands are not independent but we
can use negative association arguments (Dubhashi and Ranjan \cite
{dubhashi:ranjan:1998})
(see Section~\ref{sec:proofs}) to obtain Bennett inequalities for the
cumulated occupancy counts $K_{n,\overline{r}}$. These hold either with
the Jackknife variance proxy given by the Efron--Stein inequality,
$r\EXP K_{n,r}$ or with the variance proxy stemming from the negative
correlation of the summands, $\EXP K_{n,\overline{r}}$. Letting
$v_{n,\overline{r}}=\min(r\EXP K_{n,r}, \EXP K_{n,\overline{r}})$,
we have,
for all $\lambda\in\R$:
\[
\log\EXP\mathe^{\lambda(K_{n,\overline{r}}-\EXP K_{n,\overline
{r}})} \leq v_{n,\overline{r}} \phi(\lambda).
\]
This in turn implies a concentration inequality for $K_{n,r}$. Letting
\[
v_{n,r} = 2\min \bigl(\max\bigl(r\EXP K_{n,r}, (r+1)\EXP
K_{n,r+1}\bigr), \EXP K_{n,\overline{r}} \bigr) ,
\]
we have, for all $s\geq0$,
\[
\PROB \bigl\{ \vert K_{n,r}-\EXP K_{n,r}\vert\geq\sqrt{4
v_{n,r} s} + {2s}/3 \bigr\} \leq4 \mathrm{e}^{-s}.
\]

We obtain distribution-free bounds on the log-Laplace transform of
$M_{n,0}$, which result in sub-Gaussian concentration on the left tail,
sub-gamma concentration on the right tail with scale proxy $1/n$. More
precisely, letting $v^-_n=2\EXP K_2(n)/n^2$ and $v^+_n=2\EXP
K_{\overline{2}}(n)/n^2$, we show that, for all $\lambda\leq0$,
\[
\log\EXP\mathrm{e}^{\lambda(M_{n,0}-\EXP M_{n,0})} \leq v^-_n \frac
{\lambda^2}{2} ,
\]
and, for all $\lambda\geq0$,
\[
\log\EXP\mathrm{e}^{\lambda(M_{n,0}-\EXP M_{n,0})} \leq v^+_n \frac
{\lambda^2}{2(1-\lambda/n)}.
\]

Indeed, these results are distribution-free. But though the variance
factor $v^-_n$ is a sharp bound for the variance of the missing mass,
$v^+_n$ may be much larger. This leads us to look for
distribution-specific conditions ensuring that $v^+_n$ captures the
right order for the variance, such as by using a tail asymptotic
stability condition as in extreme value theory.

Karlin \cite{Kar67} pioneered such a condition by assuming that the function
$\vec{\nu}: (0,1]\to\mathbb{N}$, defined by $\vec{\nu}(x) =|\{j
\in
\mathbb{N}_+, p_j \geq x\}|$ satisfies a regular variation assumption,
namely $\vec{\nu}(1/n)\sim n^{\alpha} \ell(n)$ near $+\infty$, with
$\alpha\in(0,1]$ (see also Gnedin \textit{et al.} \cite
{gnedin2007notes}, Ohannessian and Dahleh \cite
{ohannessian2012rare}). Here $\ell$ is a slowly
varying function at $\infty$, i.e. for all $x$, $\ell(\tau x)/\ell
(\tau
)\to1$ as $\tau\to\infty$. This condition allows us to compare the
asymptotics of the various occupancy scores. In particular, in this
framework $\EXP K_2(n)$ and $\EXP K_{\overline{2}}(n)$ have the same
order of growth, and, divided by $n^2$ they both are of the same order
as the variance of the missing mass. Hence, regular variation provides
a framework in which our concentration inequalities are order-optimal.

To handle the case $\alpha=0$, we move from \emph{Karamata} to \emph{de
Haan} theory, and take $\vec{\nu}$ to have an extended regular
variation property, with the additional assumption that $\EXP K_1(n)$
tends to $+\infty$. This domain corresponds to light-tailed
distributions which are still heavier than the geometric. In this case,
we manage to show the sub-gamma concentration of the missing mass only
for $n$ large enough, that is, that there exists $n_0$ such that for
all $n\geq n_0$, for $\lambda>0$, we have $\log\EXP\mathrm
{e}^{\lambda(M_{n,0}-\EXP M_{n,0})} \leq(v_n\lambda^2)\slash
2(1-\lambda/n)$, with $v_n \asymp\var M_{n,0}$.

Back to our examples of applications, considerable insight may be
gained from these concentration results. For instance, heavy tails lead
to multiplicative concentration for $M_{n,0}$ and the consistency of
the Good--Turing estimator: $\frac{G_{n,0}}{M_{n,0}} \overset
{p}{\rightarrow} 1$. Generally, new estimators can be derived and
shown to be consistent in a unified framework, once one is able to
estimate $\alpha$ consistently. For instance, when $\vec\nu(1/\cdot)$ is
regularly varying with index $\alpha$, $\hat\alpha=K_{n,1}/K_n$ is a
consistent estimator of $\alpha$. Then, to estimate the number of new
species in a sample twice the size of the original, we immediately get
that $\widehat{K}_{2n}=K_n+\frac{2^{\widehat\alpha}-1}{\widehat
\alpha}
K_{n,1}$ is a consistent estimator of $K_{2n}$. This methodology is
very similar to extreme value theory (Beirlant \textit{et al.} \cite
{BeGoSeTe06}): harnessing
limiting expressions and tail parameter estimation. These results
strengthen and extend the contribution of Ohannessian and Dahleh \cite
{ohannessian2012rare},
which is restricted to power-laws and implicit constants in the
inequalities. Beyond consistency results, we also obtain confidence
intervals for the Good--Turing estimator in the Poisson setting, using
empirical quantities.

\textit{Historical notes and related work}.
There exists a vast literature on the occupancy scheme, as formulated
here and in many other variations. The most studied problems are the
asymptotic behavior of $K_n$ and $K_{n,r}$. This is done often in a
finite context, or a scaling model where probabilities remain mostly
uniform. Of particular relevance to this paper, we mention the work of
Karlin~\cite{Kar67}, who built on earlier work by Bahadur \cite{Bahadur1960},
credited as one of the first to study the infinite occupancy scheme.
Karlin's main results were to establish central limit theorems in an
infinite setting, under a condition of regular variation. He also
derived strong laws of large numbers. Gnedin \textit{et al.} \cite
{gnedin2007notes} present a
general review of these earlier results, as well as more contemporary
work on this problem. The focus continues to be central limit theorems,
or generally asymptotic results. For example, the work of Hwang and
Janson \cite
{MR2408581} (effectively) provides a local limit theorem for $K_n$
provided that the variance tends to infinity. Somewhat less asymptotic
results have also been proposed, in the form of normal approximations,
such as in the work of Bogachev \textit{et al.} \cite{MR2412154} and Barbour and Gnedin \cite{MR2480545}.

Besides occupancy counts analysis, a distinct literature investigates
the number of species and missing mass problems. These originated in
the works of Fisher \textit{et al.} \cite{Fisher1943}, Good \cite{Good1953}, and Good and
Toulmin \cite
{GoodToulmin1956}, and generated a long line of research to this day
(Bunge and Fitzpatrick \cite{bunge1993estimating}). Here, instead of
characterizing the
asymptotic behavior of these quantities, one is interested in
estimating $K_{\lambda n}-K_n$ for a $\lambda>1$, that is the number of
discoveries when the sample size is multiplied by $\lambda$, or
estimating $M_{n,0}$: estimators are proposed, and then their
statistical properties are studied. One recently studied property by
McAllester and Schapire \cite{McAllesterSchapire2000} and {McAllester}
and Ortiz \cite{mcallester2003concentration},
is that of concentration, which sets itself apart from the CLT-type
results in that it is non-asymptotic in nature. Based on this,
Ohannessian and Dahleh \cite
{ohannessian2012rare} showed that in the regular variation setting of
Karlin, one could show multiplicative concentration, and establish
strong consistency results. Conversely, characterizing various aspects
of concentration allows one to systematically design new estimators.
For example, this was illustrated in Ohannessian and Dahleh \cite
{ohannessian2012rare} for
the estimation of rare probabilities, to both justify and extend Good's
(Good \cite{Good1953}) work that remains relevant in some of the
aforementioned applications, especially computational linguistics.
Such concentration results for rare probabilities have been also used
in the general probability estimation problem, such as by Acharya \textit{et al.}
\cite{Ach}.

\section{Distribution-free concentration}
\label{sec:concentration}

\subsection{Occupancy counts}
\label{sec:conc-ineq-occup}

\subsubsection{Variance bounds}
\label{sec:vari-bounds-occup}

In order to understand the fluctuations of occupancy counts $K_n$,
$K(t)$, $K_{n,r}$, $K_r(t)$, we start by reviewing and stating variance
bounds. We start with the Poisson setting where occupancy counts are
sums of independent Bernoulli random variables with possibly different
success probabilities, and thus variance computations are
straightforward. Exact expressions may be derived (see, for
example, Gnedin \textit{et al.} \cite{gnedin2007notes}, equation
(4)), but ingenuity may be used to
derive more tractable and tight bounds. We start by stating a
well-known connection between the variance of the number of occupied
urns and the expected number of singletons (Gnedin \textit{et al.}
\cite{gnedin2007notes},
Karlin \cite{Kar67}).
In the binomial setting, similar bounds can be derived using the
Efron--Stein--Steele inequalities, outlined in Section~\ref
{sec:efron-stein} (see Boucheron \textit{et al.} \cite{BoLuMa13},
Section~3.1).

\begin{prop}
\label{prop:varKt}
In the Poissonized setting, we have
\[
\frac{\EXP K_1(2t) }{2} \leq\var\bigl(K(t)\bigr) \leq\EXP K_1(t).
\]
%
In the binomial setting, we have
%
\[
\var(K_n) \leq\EXP \bigl[ K_{n,1} (1-M_{n,0})
\bigr] \leq\EXP K_{n,1}.
\]
\end{prop}

The upper bounds in these two propositions parallel each other, but in
the binomial setting, we cannot hope to establish lower bounds like
$\EXP K_{cn,1}/c\leq\var(K_n)$ for some constant $c>0$ in full
generality. To see this, consider the following example which shows
that the variance of $K(t)$ and of $K_n$ may differ significantly, and
that the variance of $K_n$ may be much smaller than the expected value
of $K_{n,1}$.

\begin{ex} \label{ex:birthday1}
In the so-called \emph{birthday paradox} scenario, $n$ balls are thrown
independently into $n^2$ urns with uniform probabilities $1/n^2$. In
the Poisson setting for time $t=n$, since we have $\EXP K(n)=\sum_j
(1-\mathe^{-np_j}) = n^2(1-\mathe^{-1/n})$, using an upper and lower
Taylor expansion we can obtain the bounds:
\[
n-\frac{1}{2}\leq\EXP K(n) \leq n-\frac{1}{2}+\frac{1}{6n}.
\]
Since $\var(K(n))=\sum_j \mathe^{-np_j}(1-\mathe^{-np_j}) = \EXP
K(2n)-\EXP K(n)$, it follows that:
\[
n-\frac{1}{6n} \leq\var\bigl(K(n)\bigr) \leq n+\frac{1}{12 n}.
\]
Meanwhile, we have $\EXP[K_1(n)]=\sum_j np_j \mathe^{-np_j}= n\mathe
^{-1/n}$, which can be bounded similarly:
\[
n-1 \leq\EXP\bigl[K_1(n)\bigr]\leq n-1+\frac{1}{2n}.
\]
We can thus see that the Poisson birthday paradox satisfies the spirit
of Proposition~\ref{prop:varKt}, if not its letter (because of being
outside of our fixed-$p$ setting). Namely, the Poisson quantities $\var
(K(n))$ and $\EXP K_1(n)$ are of the same order of magnitude, roughly $n$.

On the other hand, in the binomial setting, since $1- M_{n,0}=K_n /n^2$,
\[
\var(K_n)\leq\EXP \bigl[ K_{n,1} (1-M_{n,0})
\bigr]=\EXP \biggl[K_{n,1}\frac{K_n}{n^2} \biggr] \leq
\frac{1}{n} \EXP K_{n,1} \leq 1, 
\]
where we have used the same variance bound as in the proof of
Proposition~\ref{prop:varKt} (Section~\ref
{sec:proof-vari-bounds-occ-counts}). While this implies that the upper
bound $\var(K_n)\leq\EXP K_{n,1}$ is satisfied, it also shows that a
lower bound of the kind of $\EXP K_{cn,1}/c\leq\var(K_n)$ is not
possible, since similarly to $\EXP K_1(n)$, we have $\EXP K_{n,1}=\sum_j n p_j (1-p_j)^n = n(1-\frac{1}{n^2})^n \geq n-1$.
\end{ex}

Another straightforward bound on $\var(K_n)$ comes from the fact that
the Bernoulli variables $(\IND_{\{X_{n,j}>0\}})_{j\geq1}$ are
negatively correlated. Thus, ignoring the covariance terms, we get
\[
\var(K_n) \leq\sum_{j=1}^\infty
\var(\IND_{\{X_{n,j}>0\}}) = \sum_{j=1}^\infty(1-p_j)^n
\bigl(1-(1-p_j)^n\bigr) =\EXP K_{2n} - \EXP
K_n.
\]
Let us denote this bound by $\operatorname{Var}^{\mathrm
{ind}}(K_n)=\EXP K_{2n} - \EXP K_n$, as
it is a variance proxy obtained by considering that the summands in
$K_n$ are independent. One can observe that the expression of
$\operatorname{Var}^{\mathrm{ind}}
(K_n)$ is very similar to the variance in the Poissonized setting,
$\var
(K(t))=\EXP K(2t)-\EXP K(t)$. It is insightful to compare the true
variance, the Poissonized proxy, and the negative correlation proxy, to
quantify the price one pays by resorting to the latter two as an
approximation for the first. We revisit this in more detail in
Section~\ref{sec:cost-poisson}.

We now investigate the fluctuations of the individual occupancy counts
$K_{n,r}$ and $K_r(t)$, and that of the cumulative occupancy counts
$K_{n,\overline{r}}= \sum_{j\geq r} K_{n,j}$ and $K_{\overline{r}}(t)=
\sum_{j\geq r} K_j(t)$.
%
\begin{prop}
\label{prop:varKnrinf}
In the Poisson setting, for $r\geq1$, $t\geq0$,
\[
\var\bigl(K_{\overline{r}}(t)\bigr) \leq\min \bigl(r \EXP K_r(t),
\EXP K_{\overline{r}}(t) \bigr).
\]
In the binomial setting, for $r,n\geq1$,
\[
\var(K_{n,\overline{r}}) \leq\min (r \EXP K_{n,r}, \EXP
K_{n,\overline{r}} ) .
\]
\end{prop}
For each setting, the first bound follows from Efron--Stein--Steele
inequalities, the second from negative correlation. These techniques
are presented briefly in Sections~\ref{sec:efron-stein} and \ref
{sec:negative-association}, respectively.

\begin{rem}
Except for $r=1$, there is no clear-cut answer as to which of these
two bounds is the tightest. In the regular variation scenario with
index $\alpha\in\,]0,1]$ as explored in Gnedin \textit{et al.}~\cite
{gnedin2007notes}, the
two bounds are asymptotically of the same order, indeed,
\[
\frac{r\EXP K_{n,r}}{\EXP K_{n,\overline{r}}}\mathop{\sim}_{n\to+\infty} \alpha,
\]
see Section~\ref{sec:regular-variation} for more on this.
\end{rem}

Bounds on $\var(K_r(t))$ can be easily derived as $K_{r}(t)$ is a sum
of independent Bernoulli random variables, however, noticing that
$K_{n,r}=K_{n,\overline{r}}-K_{n,\overline{r+1}}$ and that
$K_{n,\overline{r}}$ and $K_{n,\overline{r+1}}$ are positively
correlated, the following bound is immediate.
%
\begin{prop}
\label{prop:varKnr}
In the Poisson setting, for $r\geq1$, $t\geq0$,
\[
\var\bigl(K_r(t)\bigr) \leq\EXP K_r(t).
\]
In the binomial setting, for $r,n\geq1$,
\begin{eqnarray*}
\var(K_{n,r}) &\leq& \min \bigl( r\EXP K_{n,r} + (r+1)\EXP
K_{n,r+1}, \EXP K_{n,\overline{r}} + \EXP K_{n,\overline{r+1}} \bigr)
\\
& \leq& 2 \min \bigl(\max\bigl(r\EXP K_{n,r}, (r+1)\EXP
K_{n,r+1}\bigr), \EXP K_{n,\overline{r}} \bigr).
\end{eqnarray*}
\end{prop}

\subsubsection{Concentration inequalities}

Concentration inequalities refine variance bounds. These bounds on the
logarithmic moment generating functions are indeed Bennett
(sub-Poisson) inequalities with the variance upper bounds stated in the
preceding section. 
For $K_{n,\overline{r}}$, the next proposition gives a Bernstein inequality
where the variance factor is, up to a constant factor, the Efron--Stein
upper bound on the variance.
%
\begin{prop}
\label{prop:conc-ineq-knrinf}
Let $r\geq1$, and let $v_{n,\overline{r}}=\min(r\EXP K_{n,r}, \EXP
K_{n,\overline{r}})$. Then, for all $\lambda\in\mathbb{R}$,
\[
\log\EXP\mathe^{\lambda(K_{n,\overline{r}}-\EXP K_{n,\overline
{r}})} \leq v_{n,\overline{r}} \phi(\lambda) ,
\]
with $\phi:\lambda\mapsto\mathe^\lambda-\lambda-1$.
\end{prop}
It is worth noting that the variance bound $\EXP K_{n,\overline{r}}$
in this
concentration inequality can also be obtained using a variant of
Stein's method known as size-biased coupling (Bartroff \textit{et al.}
\cite
{bartroff2014bounded},
Chen \textit{et al.} \cite{chen2010normal}).

A critical element of the proof of Proposition~\ref
{prop:conc-ineq-knrinf} is to use the fact that each $K_{n,\overline
{r}}$ is
a sum of negatively associated random variables (Section~\ref
{sec:negative-association}). This is not the case for $K_{n,r}$, and
thus negative association cannot be invoked directly. To deal with
this, we simply use the observation of Ohannessian and Dahleh \cite
{ohannessian2012rare} that
since $K_{n,r}= K_{n,\overline{r}}-K_{n,\overline{r+1}}$, the
concentration of
$K_{n,r}$ follows from that of those two terms. We can show the following.

\begin{prop}
\label{prop:conc-ineq-knr}
Let
\[
v_{n,r} = 2\min \bigl(\max\bigl(r\EXP K_{n,r}, (r+1)\EXP
K_{n,r+1}\bigr), \EXP K_{n,\overline{r}} \bigr).
\]
Then, for $s\geq0$,
\[
\PROB \bigl\{ \vert K_{n,r}-\EXP K_{n,r}\vert\geq\sqrt{4
v_{n,r} s} + {2s}/3 \bigr\} \leq4 \mathrm{e}^{-s}.
\]
\end{prop}

\subsection{Missing mass}
\label{sec:conc-ineq-miss}

\subsubsection{Variance bound}

Recall that $M_{n,0}=\sum_{j=1}^\infty p_j \IND_{\{X_{n,j}=0\}}=\sum_{j=1}^\infty p_j Y_j$, and we can readily show that the summands are
negatively associated weighted Bernoulli random variables (Section~\ref
{sec:negative-association}). This results in a handy upper bound for
the variance of the missing mass.

\begin{prop} \label{prop:vari-missingmass}
In the Poisson setting,
\[
\var\bigl(M_0(t)\bigr)= 2\EXP K_2(t)/t^2 -
\EXP K_2(2t)/2t^2 \leq2\EXP K_2(t)/t^2
,
\]
while in the binomial setting,
\[
\var(M_{n,0}) \leq\sum_{j=1}^\infty
p_j^2 \var(Y_j) \leq\frac{2
\EXP
K_{2}(n) }{n^2}.
\]
\end{prop}
Note that whereas the expected value of the missing mass is connected
to the number of singletons, its variance may be upper bounded using
the number of \emph{doubletons} (in the Poisson setting). This
connection was already pointed out in Good \cite{Good1953} and Esty
\cite{Est82}.

\subsubsection{Concentration of the left tail}

Moving on to the concentration properties of the missing mass, we first
note that as a sum of weighted sub-Poisson random variables (following
Boucheron \textit{et al.} \cite{BoLuMa13}), the missing mass is
itself a
sub-gamma random
variable on both tails. It should not come as a surprise that the left
tail of $M_{n,0}$ is sub-Gaussian with the variance factor derived from
negative association. This had already been pointed out by McAllester
and Schapire \cite
{McAllesterSchapire2000} and {McAllester} and Ortiz \cite
{mcallester2003concentration}.

\begin{prop}[(McAllester and Ortiz \cite{mcallester2003concentration})]
\label{prop:bennett:left:mn0}
In the Poisson setting, the missing mass $M_0(t)$ is sub-Gaussian on
the left tail with variance factor the effective variance $\var
(M_0(t))=\sum_{j=1}^\infty p_j^2 \mathe^{-tp_j}(1-\mathe^{-tp_j})$.

In the binomial setting, the missing mass $M_{n,0}$ is sub-Gaussian on
the left tail with variance factor $v=\sum_{j=1}^\infty p_j^2 \var
(Y_j)$ or $v_n^-=2\EXP K_2(n)/n^2$.

For $\lambda\leq0$,
\[
\log\EXP \bigl[ \mathe^{\lambda(M_{n,0}-\EXP M_{n,0})} \bigr] \leq \frac{v \lambda^2}{2} \leq
\frac{2 v_n^- \lambda^2}{2}.
\]
\end{prop}

\subsubsection{Concentration of the right tail}

The following concentration inequalities for the right tail of the
missing mass mostly rely on the following proposition, which bounds the
log-Laplace transform of the missing mass in terms of the sequence of
expected occupancy counts $\EXP K_r(n)$, for $r\geq2$.

\begin{prop} \label{prop:log-laplace-missing-mass}
For all $\lambda>0$,
\[
\log\EXP \bigl[ \mathe^{\lambda(M_{n,0}-\EXP M_{n,0})} \bigr] \leq \sum
_{r=2}^\infty \biggl(\frac{\lambda}{n}
\biggr)^r \EXP K_r(n).
\]
\end{prop}

This suggests that if we have a uniform control on the expected
occupancy scores $(\EXP K_r(t))_{r\geq2}$, then the missing mass has a
sub-gamma right tail, with some more or less accurate variance proxy,
and scale factor $1/n$.

The next theorem shows that the missing mass is sub-gamma on the right
tail with variance proxy $2\EXP K_{\overline{2}}(n)/n^2$ and scale
proxy $1/n$. Despite its simplicity and its generality, this bound
exhibits an intuitively correct scale factor: if there exist symbols
with probability of order $1/n$, they offer the major contribution to
the fluctuations of the missing mass.

\begin{thmm}\label{thmm:struct-ineq-right}
In the binomial as well as in the Poisson setting, the missing mass is
sub-gamma on the right tail with variance factor $v_n^+= 2\EXP
K_{\overline{2}}(n)/n^2 $ and scale factor $1/n$.
For $\lambda\geq0$,
\[
\log\EXP \bigl[ \mathe^{\lambda(M_{n,0}-\EXP M_{n,0})} \bigr] \leq \frac{v_n^+\lambda^2}{2(1-\lambda/n)}
.
\]
If the sequence $(\EXP K_r(n))_{r\geq2}$ is non-increasing, the
missing mass is sub-gamma on the right tail with variance factor
$v_n^-= 2\EXP K_2(n)/n^2$ and scale factor $1/n$,
\[
\log\EXP \bigl[ \mathe^{\lambda(M_{n,0}-\EXP M_{n,0})} \bigr] \leq \frac{v_n^-\lambda^2}{2 (1-\lambda/n)}.
\]
\end{thmm}


\begin{rem}
{McAllester} and Ortiz \cite{mcallester2003concentration} and Berend
and Kontorovich \cite{BerKon13} point out that
for each Bernoulli random variable $Y_j$, for all $\lambda\in\mathbb{R}$
\[
\log\EXP\mathe^{\lambda(Y_j-\EXP Y_j)} \leq\frac{\lambda
^2}{4\textsc{c}_{\textsc{ls}}(\EXP Y_j)} ,
\]
where
$\textsc{c}_{\textsc{ls}}(p) = \log((1-p)/p)/(1-2p)$ (or $2$ if $p=1/2$)
is the optimal logarithmic Sobolev constant for Bernoulli random
variables with success probability $p$ (this sharp and nontrivial
result has been proven independently by a number of people: Kearns and
Saul \cite
{kearns1998large}, Berend and Kontorovich \cite{BerKon13}, Raginsky
and Sason \cite{RagSas12}, Berend and Kontorovich \cite
{BerKonJMLR}; the constant also appears early on in the exponent of one
of Hoeffding's inequalities (Hoeffding \cite{Hoe63}, Theorem~1, equation
(2.2)). From this observation, thanks to the negative
association of the $(Y_j)_{j\geq1}$, it follows that the missing mass
is sub-Gaussian with variance factor
%
\begin{equation}
\label{eq:7} w_n= \sum_{j=1}^\infty
\frac{p_j^2}{2\textsc{c}_{\textsc
{ls}}((1-p_j)^n) } \leq\sum_{j=1}^\infty
\frac{p_j^2}{2\log((1-p_j)^{-n})} \leq\sum_{j=1}^\infty
\frac{p_j^2}{2 np_j} \leq\frac{1}{2n}.
\end{equation}
An upper bound on $w_n$ does not mean that $w_n$ is necessarily larger
than $\EXP K_{\overline{2}}(n)/n^2$. Nevertheless, it is possible to
derive a simple lower bound on $w_n$ that proves to be of order
$O(1/n)$.

Assume that the sequence $(p_j)_{j\geq1}$ is such that
$p_j \leq1/4$ for all $j\geq1$. Then
\begin{eqnarray*}
w_n & \geq& \sum_{j: p_j\geq1/n}
\frac{p_j^2}{2 \textsc
{c}_{\textsc
{ls}}((1-p_j)^n)}
\\
& \geq& \sum_{j : p_j\geq1/n} \frac{p_j^2 (1-2(1-p_j)^n)}{2 n \log
(1/(1-p_j))}
\\
& \geq& \sum_{j : p_j\geq1/n} \frac{p_j^2 (1-2/\mathrm{e})}{2 n
p_j/(1-p_j)}
\\
& \geq& \frac{3(1-2/\mathrm{e})}{8 n}\biggl(1- \sum_{j : p_j<1/n}
p_j\biggr)
\\
& \geq& \frac{3}{32n} \biggl(1- \sum_{j : p_j<1/n}
p_j \biggr) ,
\end{eqnarray*}
and the statement follows from the observation that $\lim_{n\to\infty}
\sum_{j : p_j<1/n} p_j =0$.

The variance factor $w_n$ from \eqref{eq:7} is usually larger than
$2\EXP
K_{\overline{2}}(n)/n^2$. In the scenarios discussed in Section~\ref
{sec:regular-variation}, $(2\EXP K_{\overline{2}}(n)/n^2)/w_n$
even tends to $0$ as $n$ tends to infinity.
\end{rem}

\section{Regular variation}
\label{sec:regular-variation}

\textit{Motivation}.
Are the variance bounds in the results of Section~\ref
{sec:concentration} tight? In some pathological situations, this may
not be the case. Consider the following example (which revisits Example~\ref{ex:birthday1}).

\begin{ex} \label{ex:birthday2}
We may challenge the tail bounds offered by Proposition~\ref
{prop:bennett:left:mn0} and Theorem~\ref{thmm:struct-ineq-right} in the
simplest setting where we have $k$ symbols all of which have equal
probabilities $1/k$. Then the missing mass is $1-K_n/k$, its variance
is $\var(K_n) /k^2$. In the birthday paradox setting ($k=n^2$), $\var
(M_{n,0})\leq1/n^4$, and the variance bound $2 \EXP K_{\overline
{2}}(n)/n^2$ is not tight. Indeed, one can verify that $\EXP
K_{\overline{2}}(n)\geq\frac{1}{2} (1-\frac{1}{n} )$ so
that $v_n^+\geq\frac{1}{n^2}-\frac{1}{n^3}$. However, in what is
called the \emph{central domain} in Kolchin \textit{et al.} \cite{MR0471016},
that is when $k
\to\infty$ while $n/k\to t \in\mathbb{R}_+$, the tail bounds become
relevant. The variance of $K_n$ is equivalent to $k \mathe^{-t}
(1-\mathe^{-t})$ while its expectation is equivalent to $k (1-\mathe
^{-t})$. Note that in this setting all $\EXP K_r(n)$ and $\EXP K_{n,r}$
are of the same order of magnitude as $\EXP K_n$, indeed $\EXP
K_r(n)/\EXP K(n) \to\mathe^{-t} t^r/(r! (1-\mathe^{-t}))$.

These examples are illustrative although they do not fall in the fixed
$p$ regime we are considering in this paper. We use them because they
have tractable expressions, and they provide informative diagnostics.
To parallel the phenomenon of mismatched variance proxies in our
setting, one can simply look at the geometric distribution for a
concrete example. If $(p_k)_{k\geq1}$ defines a geometric distribution
$p_k=(1-q)^{k-1}q$, then $\EXP K_2(n)$
remains bounded, while $\EXP K_{\overline{2}}(n)$ scales like $\log n$
as $n$ tends to infinity.
\end{ex}

In particular, we may conjecture that Theorem~\ref
{thmm:struct-ineq-right} is likely to be sharp when the first terms of
the sequence $(\EXP K_r(n))_{r\geq2}$ grow at the same rate as $\EXP
K(n)$, or at least as $\EXP K_{\overline{2}}(n)$, which is not the case
in the birthday paradox setting of Example~\ref{ex:birthday2}. We see
in what follows that the regular variation framework introduced by
Karlin \cite{Kar67} leads to such asymptotic equivalents. The most useful
aspect of these equivalent growth rates is a simple characterization of
the variance of various quantities, particularly relative to their
expectation. We focus on the right tail of the missing mass, which
exhibits the highest sensitivity to this asymptotic behavior, by trying
to specialize Theorem~\ref{thmm:struct-ineq-right} under regular variation.

\textit{Definition}.
Regularly varying frequencies can be seen as generalizations of
power-law frequencies. One possible definition is as follows: for
$\alpha\in(0,1)$, the sequence $(p_{j})_{j\geq1}$ is said to be
regularly varying with index $-1/\alpha$ if, for all $\kappa\in\N_{+}$,
\[
\frac{p_{\kappa j}}{p_{j}}\mathop{\sim}\limits_{j\to\infty}\kappa^{-{1}/{\alpha}} .
\]
It is easy to see that pure power laws do indeed satisfy this
definition. However, in order to extend the regular variation
hypothesis to $\alpha=0$ and $1$, we need a more flexible definition,
which requires some new notation.
This definition relies on the \textit{counting function} $\vec\nu$,
defined for all $x>0$ by:
\[
\vec{\nu}(x) =\bigl\llvert \{ j : j \geq1 , p_j\geq x \} \bigr
\rrvert .
\]
The overhead arrow is not a vector notation, and rather codifies that
we are counting points with probability ``to the right'' of $x$. More
precisely, letting $\nu$ be the \textit{counting measure} defined by
\[
\nu(\mathrm{d}x)=\sum_{j=1}^\infty
\delta_{p_j}(\mathrm{d}x) ,
\]
then, for all $x>0$, $\vec\nu(x)=\nu[x,1]$.

Henceforth, following Karlin \cite{Kar67}, we say that the probability mass
function $(p_j)_j$ is regularly varying with index $\alpha\in
[0,1]$, if $\vec{\nu}(1/\cdot)$ is $\alpha$-regularly varying in the
neighbourhood of $\infty$, which reads as
\[
\vec{\nu}(1/x)\mathop{\sim}_{x \to\infty} x^{\alpha}\ell(x) ,
\]
where $\ell$ is a slowly varying function, that is, for all $x>0$,
$\lim_{\tau\to+\infty} {\ell(\tau x)}/{\ell(\tau)}=1$. We use the notation
$\vec\nu(1/\cdot)\in\textsc{RV}_\alpha$ to indicate that $\vec
{\nu
}(1/\cdot)$ is $\alpha$-regularly varying.

We now note that when $\alpha\in(0,1)$, the regular variation
assumption on $(p_j)_{j\geq1}$ is indeed equivalent to the regular
variation assumption on the counting function $\vec{\nu}$ (see Gnedin et
al.~\cite{gnedin2007notes}, Proposition~23): if $(p_j)_{j\geq1}$ is
regularly varying with index $-1/\alpha$ as $j$ tends to infinity, then
$\vec{\nu}(1/\cdot)$ is $\alpha$-regularly varying, that is $\lim_{x\to
\infty} \vec{\nu}(1/(qx))/\vec{\nu}(1/x)=q^{\alpha}$ for $q>0$
(see also Bingham \textit{et al.} \cite{BiGoTe87}). The second
definition however lends itself more
easily to generalization to $\alpha=0$ and $1$.

In what follows, we treat these three cases separately: the nominal
regular variation case with $\alpha\in(0,1)$ strictly, the \emph{fast
variation} case with $\alpha=1$, and the \emph{slow variation} case
with $\alpha=0$.

In the latter case, that is if frequencies $p_j$ are regularly varying
with index $0$, we find that the mere regular variation hypothesis is
not sufficient to obtain asymptotic formulas. For this reason, we
introduce further control in the form of an \emph{extended} regular
variation hypothesis (given by Definition~\ref{dfn:slow-variation-1} of
Section~\ref{sec:slow-variation}). 

\begin{rem} \label{rem:necessary}
Before we proceed, as further motivation, we note that the regular
variation hypothesis is very close to being a necessary condition for
exponential concentration. For example, considering Proposition~\ref
{prop:log-laplace-missing-mass}, we see that if the sampling
distribution is such that the ratio $\EXP K_{\overline{2}}(t)/ \EXP
K_2(t)$ remains bounded, then we are able to capture the right variance
factor. Now, defining the shorthand $\Phi_{\overline{2}}(t)=\EXP
K_{\overline{2}}(t)$ and $\Phi_2(t)=\EXP K_2(t)$ following the notation
of Gnedin \textit{et al.} \cite{gnedin2007notes}, we have
\[
\Phi_{\overline{2}}'(t)=\frac{2\Phi_2(t)}{t}.
\]
Hence, $\Phi_{\overline{2}}(t)/ \Phi_2(t)=2\Phi_{\overline
{2}}(t)/t\Phi
_{\overline{2}}'(t)$, and if instead of boundedness, we further require
that this ratio converges to some finite limit, then, by the converse
part of Karamata's theorem (see {de Haan} and Ferreira \cite
{HaaFei06}, Theorem
B.1.5), we
find that $\Phi_{\overline{2}}$ (and then $\Phi_2$) is regularly
varying, which in turn implies that $\vec\nu(1/t)$ is regularly
varying. We elaborate on this further in our discussions, in
Section~\ref{sec:extensions}. 
\end{rem}

\subsection{Case \texorpdfstring{$\alpha\in(0,1)$}{$alpha in(0,1)$}}

We first consider the case $0<\alpha< 1$. The next theorem states that
when the sampling distribution is regularly varying with index $\alpha
\in(0,1)$, the variance factors in the Bernstein inequalities of
Proposition~\ref{prop:bennett:left:mn0} and Theorem~\ref
{thmm:struct-ineq-right} are of the same order as the variance of the
missing mass.

\begin{thmm}\label{prop:bernstein:right:regvar}
Assume that the counting function $\vec{\nu}$ satisfies the regular
variation condition with index $\alpha\in(0,1)$, then the missing mass
$M_{n,0}$ (or $M_0(n)$) is sub-Gaussian on the left tail with variance
factor $v_n^-=2\EXP K_2(n)/n^2$ and sub-gamma on the right tail with
variance factor $v_n^+=2\EXP K_{\overline{2}}(n)/n^2$. The variance
factors satisfy
\begin{eqnarray*}
\lim_n \frac{v_n^-} {\var(M_{n,0})}&=& \frac{1}{1-2^{\alpha-2}} ,
\\
\lim_n \frac{v_n^+} {\var(M_{n,0})}&= &\frac{2}{\alpha(1-2^{\alpha-2})} ,
\end{eqnarray*}
and thus
\[
\lim_n \frac{v_n^-} {v_n^+}= \frac{\alpha}{2}.
\]
\end{thmm}

The second ratio deteriorates when $\alpha$ approaches $0$, implying
that the variance factor for the right tail gets worse for lighter
tails. We do not detail the proof of Theorem~\ref
{prop:bernstein:right:regvar}, except to note that it follows from
Proposition~\ref{prop:bennett:left:mn0}, Theorem~\ref
{thmm:struct-ineq-right}, and the following asymptotics (see
also Gnedin \textit{et al.}
\cite{gnedin2007notes}, Ohannessian and Dahleh \cite{ohannessian2012rare}).
%
\begin{thmm}[(Karlin \cite{Kar67})]
\label{th:portmanteau}
If the counting function $\vec{\nu}$ is regularly varying with index
$\alpha\in(0,1)$, for all $r\geq1$,
\begin{itemize}[--]
\item[--]$K_n \overset{a.s.}\sim\EXP K_n \sim_{+\infty}\Gamma
(1-\alpha)n^\alpha\ell(n) $,
\item[--]$K_{n,r} \overset{a.s.}\sim\EXP K_{n,r}\sim_{+\infty}
\frac{\alpha\Gamma(r-\alpha)}{r!}n^\alpha\ell(n) $,
\item[--]$\var(M_{n,0}) \sim\alpha\Gamma(2-\alpha)(1-2^{\alpha
-2})n^{\alpha-2}\ell(n) $
\end{itemize}
and the same hold for the corresponding Poissonized quantities.
\end{thmm}
Note that all expected occupancy counts are of the same order, and the
asymptotics for $\EXP K_{\overline{2}}(n)$ follows directly from the
difference between $\EXP K(n)$ and $\EXP K_1(n)$.

\subsection{Fast variation, \texorpdfstring{$\alpha=1$}{$alpha=1$}}
\label{sec:fast-variation}

We refer to the regular variation regime with $\alpha=1$ as \emph{fast
variation.}\footnote{Sometimes \emph{rapid variation} is used Gnedin
\textit{et al.} \cite
{gnedin2007notes}, but this conflicts with Bingham \textit{et al.}
\cite
{BiGoTe87}.} From the
perspective of concentration, this represents a relatively ``easy''
scenario. In a nutshell, this is because the variance of various
quantities grows much slower than their expectation.

The result of this section is to simply state that Theorem~\ref
{prop:bernstein:right:regvar} continues to hold as is for $\alpha=1$.
The justification for this, however, is different. In particular, the
asymptotics of Theorem~\ref{th:portmanteau} do not apply: the number of
distinct symbols $K_n$ and the singletons $K_{n,1}$ continue to have
comparable growth order, but now their growth dominates that of
$K_{n,r}$ for all $r\geq2$. Intuitively, under fast variation almost
all symbols appear only once in the observation, with only a vanishing
fraction of symbols appearing more than once. We formalize this in the
following theorem.


\begin{thmm}[(Karlin \cite{Kar67})]
\label{thmm:portmanteau-rapid-var}
Assume $\vec{\nu}(1/x)= x\ell(x)$ with $\ell\in\textsc{RV}_0$ (note
that $\ell$ tends to $0$ at $\infty$).
Define $\ell_1: [1,\infty) \to\mathbb{R}_+$ by
$
\ell_1(y)=\int_y^\infty u^{-1}\ell(u)\,\mathrm{d}u $.
Then $\ell_1 \in\textsc{rv}_0$ and $\lim_{t\to\infty} \ell
_1(t)/\ell
(t) =\infty$
and the following asymptotics hold:
\begin{itemize}[--]
\item[--]$K_n \overset{a.s.}\sim\EXP K_n \sim_{+\infty} n\ell_1(n)$,
\item[--]$K_{n,1} \overset{a.s.}\sim\EXP K_{n,1} \sim_{+\infty}
\EXP K_n$,
\item[--]$K_{n,r} \overset{a.s.}\sim\EXP K_{n,r}\sim_{+\infty}
\frac{1}{r(r-1)}n\ell(n)$, $r\geq2$,
\end{itemize}
and the same hold for the corresponding Poissonized quantities.
\end{thmm}

As the expected missing mass scales like $\EXP K_1(n)/n$ while its
variance scales like $\EXP K_2(n)/n^2$, Theorem~\ref
{thmm:portmanteau-rapid-var} quantifies our claim that this is an
``easy'' concentration. To establish Theorem~\ref
{prop:bernstein:right:regvar}, it remains to show that $\EXP
K_{\overline{2}}(n)$ is also of the same order as $\EXP K_2(n)$, with
the correct limiting ratio for $\alpha=1$. For this, we give the
following proposition, which is in fact sufficient to prove Theorem~\ref
{prop:bernstein:right:regvar} for both $0<\alpha<1$ and $\alpha=1$.


\begin{prop}\label{prop:asymp-cumulated-occup}
Assume that the counting function $\vec{\nu}$ satisfies the regular
variation condition with index $\alpha\in(0,1]$, then for all $r\geq2$,
\[
K_{\overline{r}}(n) \mathop{\sim}\limits_{+\infty} \frac{\Gamma(r-\alpha
)}{(r-1)!} \vec\nu(1/n)\qquad
\mbox{almost surely.}
\]
\end{prop}

Thus, when $\alpha=1$, $\EXP K_r(n)$ and $\EXP K_{\overline{r}}(n)$ for
$r\geq2$ all grow like $n\ell(n)$, which is dominated by the $n\ell
_1(n)$ growth of $\EXP K(n)$ and $\EXP K_1(n)$, as $\ell(n)/\ell_1(n)
\to0$. Specializing for $r=2$, we do find that our proxies still
capture the right order of the variance of the missing mass, and that
we have the desired limit of Theorem~\ref{prop:bernstein:right:regvar},
$\lim_n v_n^-/v_n^+ = \frac{1}{2}$.

\begin{rem}
When $0<\alpha<1$, another good variance proxy would have been
$2\EXP K(n)/n^2$. For $\alpha=1$, however, singletons should be removed
to get the correct order.

We also note that when $\alpha=1$, the missing mass is even more
stable. If we let $v_n$ denote either $2\EXP K_2(n)/n^2$ or $2\EXP
K_{\overline{2}}(n)/n^2$, then we have the following comparison between
the expectation and the fluctuations of the missing mass, with the
appropriate constants:
\[
\frac{\sqrt{v_n}}{\EXP M_{n,0}}\sim %
\cases{\displaystyle c_\alpha\cdot n^{-\alpha/2}
\frac{\sqrt{\ell(n)}}{\ell(n)}, &\quad $\mbox{for $0<\alpha<1$}$ , \vspace*{2pt}
\cr
\displaystyle c_1
\cdot n^{-1/2}\frac{\sqrt{\ell(n)}}{\ell_1(n)} ,&\quad $\mbox {for $\alpha=1$} $.}
%
\]
\end{rem}

\subsection{Slow variation, \texorpdfstring{$\alpha=0$}{$alpha=0$}}
\label{sec:slow-variation}

The setting where the counting function $\vec{\nu}$ satisfies the
regular variation condition with index $0$ represents a challenge. We
refer to this regime simply as \emph{slow variation}. Recall that this
means that $\vec\nu(z/n)/\vec\nu(1/n)$ converges to $1$ as $n$ goes to
infinity, yet to deal with this case we need to control the speed of
this convergence, exemplified by the notion of extended regular
variation that was introduced by de Haan (see Bingham \textit{et al.}
\cite{BiGoTe87},
de Haan and Ferreira \cite{HaaFei06}). As we illustrate in the end
of this
section, one may face rather irregular behavior without such a hypothesis.

%
\begin{dfn}\label{dfn:slow-variation-1}
A measurable function $\ell: \mathbb{R}^+ \to\mathbb{R}^+$ has the
extended slow variation property, if there exists a nonnegative
measurable function
$a: \mathbb{R}^+ \to\mathbb{R}^+$ such that for all $x>0$
\[
\lim_{\tau\to\infty} \frac{\ell(\tau x)-\ell(\tau)}{a(\tau
)}\mathop{\rightarrow}_{\tau\to\infty}
\log(x).
\]
The function $a(\cdot)$ is called an auxiliary function. When a
function $\ell$ has the extended slow variation property with auxiliary
function $a$, we denote it by $\ell\in\Pi_{a}$.
\end{dfn}
Note that the auxiliary function is always slowly varying and grows
slower than the original function, namely it satisfies $\lim_{\tau\to
\infty} \ell(\tau)/a(\tau)=\infty$. Furthermore, any two possible
auxiliary functions are asymptotically equivalent, that is if $a_1$ and
$a_2$ are both auxiliary functions for $\ell$, then $\lim_{t\to
\infty}
a_1(t)/a_2(t)=1$.

The notion of extended slow variation and the auxiliary function give
us the aforementioned control needed to treat the $\alpha=0$ case on
the same footing as the $0<\alpha<1$ case. In particular, in what
follows in this section we assume that $\vec\nu(1/\cdot)\in\Pi_a$,
with the
additional requirement that the auxiliary function $a$ tends to
$+\infty$.

\begin{rem}
This domain corresponds to light-tailed distributions just above the
geometric distribution (the upper-exponential part of Gumbel's domain).
For the geometric distribution with frequencies $p_j=(1-q)q^{k-1}$,
$j=1,2,\ldots,$ the counting function satisfies $\vec\nu(1/n)\sim
_{\infty
} \log_{1/q} (n) \in\textsc{RV}_0$, but the auxiliary function
$a(n)=\log(1/q)$ does not tend to infinity. Frequencies of the form
$p_j=cq^{\sqrt{j}}$ on the other hand do fit this framework.
\end{rem}

\begin{thmm}[(Gnedin \textit{et al.} \cite{gnedin2007notes})]
\label{thmm:portmanteau-slow-var}
Assume that $\ell(t)=\vec{\nu}(1/t)$ is in $\Pi_{a}$, where $a$ is
slowly varying and tends to infinity. The
following asymptotics hold for each $r\geq n$:
\begin{itemize}[--]
\item[--]$K_n \overset{\mathbb{P}}\sim
\EXP K_n \sim_{+\infty}\ell(n)$,
\item[--]$K_{n,r}\overset{\mathbb{P}}\sim\EXP K_{n,r} \sim_{+\infty
}\frac{a(n)}{r}$,
\item[--]$K_{n,\bar{r}}\overset{\mathbb{P}}\sim\EXP K_{n,\bar{r}}
\sim_
{+\infty}\ell(n)$,
\item[--]$M_{n,r}\overset{\mathbb{P}}\sim
\EXP M_{n,r} \sim_{+\infty}\frac{a(n)}{n}$.
\end{itemize}
The same equivalents hold for the corresponding Poissonized quantities.
\end{thmm}

\begin{rem}
In this case, the expectations $(\EXP K_{n,r})_{r\geq1}$ are of the
same order but are much smaller than $\EXP K_n$, and the variables
$K_n$ and $K_{n,\overline{s}}$ are all almost surely equivalent to
$\ell
(n)$. It is also remarkable that all the expected masses $(\EXP
M_{n,r})_{r\geq1}$ are equivalent. 
\end{rem}

The variance of the missing mass is of order $2\EXP K_{2}(n)/n^2 \sim
a(n)/n^2$, whereas the proxy $2\EXP K_{\overline{2}}(n)/n^2$ is of much
faster order $2\ell(n)/n^2$, and is thus inadequate. By exploiting more
carefully the regular variation hypothesis, we obtain uniform control
over $(\EXP K_r(n))_{r\geq1}$ for large enough $n$, leading to a
variance proxy of the correct order.


\begin{thmm}\label{thmm:bernstein:right:slowvar} Assume that $\ell$
defined by $\ell(x) = \vec{\nu}(1/x)$ is in $\Pi_{a}$ where the slowly
varying function $a$ tends to infinity, and let $v_n= {12a(n)}/n^2$. We have:
\begin{longlist}[1.]
\item[1.]$\var(M_{n,0})\sim\frac{3a(n)}{4n^2}$, thus $v_n\asymp\var(M_{n,0})$.
\item[2.] There exists $n_0\in\mathbb{N}$ that depends on $\vec{\nu}$ such
that for all $n>n_0$, for all $\lambda>0$,
\[
\log\EXP \bigl[ \mathrm{e}^{\lambda(M_{n,0}-\EXP M_{n,0})} \bigr] \leq\frac{v_n \lambda^2}{2 (1-\lambda/n)}.
\]
\end{longlist}

The same results hold for $M_0(t)$.
\end{thmm}

\begin{rem}
By standard Chernoff bounding, Theorem~\ref
{thmm:bernstein:right:slowvar} implies that there exists $n_0\in
\mathbb
{N}$ such that for all $n\geq n_0$, $s\geq0$,
\[
\PROB \biggl\{ M_{n,0} \geq\EXP M_{n,0} + \sqrt{2
v_n s} + \frac{s}{n} \biggr\} \leq \mathe^{-s}.
\]
\end{rem}

\subsubsection{Too slow variation}
\label{sec:too-slow}

We conclude this section by motivating why it is crucial to have a
heavy-enough tail in order to obtain meaningful concentration. For
example, even under regular variation when $\alpha=0$, but $\vec\nu$ is
not in a de Haan class $\Pi_{a}$ with $a(n)\to\infty$, the behavior of
the occupancy counts and their moments may be quite irregular. In this
section, we collect some observations on those light-tailed
distributions. We start with the geometric distribution which
represents in many respects a borderline case.

The geometric case is an example of slow variation: $\vec{\nu
}(1/\cdot)
\in\textsc{RV}_0$. Indeed, with $p_k=(1-q)^{k-1}q$, $0<q<1$, we have
\begin{eqnarray*}
\vec{\nu}(x)&=&\sum_{k=1}^{+\infty}
\IND_{\{p_k\geq x\}}
\\
&=& \bigl|k\in\N_+, (1-q)^{k-1}q\geq x\bigr|
\\
&=& 1+ \biggl\lfloor\frac{\log(x/q)}{\log(1-q)} \biggr\rfloor,
\end{eqnarray*}
and thus $\vec{\nu}(x)\mathrel{{\sim}_{x\to0}}\ell(1/x)$, with $\ell$
slowly varying.

In this case, $\var(K(n))=\EXP K(2n)-\EXP K(n)\to\frac{\log
(2)}{\log
 (1/1-q )}$.

\begin{prop}
\label{geom}
When the sampling distribution is geometric with parameter $q\in
(0,1)$, letting
$M_n=\max(X_1,\ldots,X_n)$,
\[
\EXP M_n \geq\EXP K_n \geq\EXP M_n -
\frac{1-q}{q^2}.
\]
\end{prop}

In the case of geometric frequencies, the missing mass can fluctuate
widely with respect to its expectation, and one cannot expect to obtain
sub-gamma concentration with both the correct variance proxy and scale
factor $1/n$. Indeed, intuitively, the symbol which primarily
contributes to the missing mass' fluctuations, is the quantile of order
$1- 1/n$. With $F(k)=\sum_{j=1}^k p_j$, and $F^{\leftarrow}$ the
generalized inverse of $F$,
\begin{eqnarray*}
j^*=F^{\leftarrow}(1-1/n)&=&\inf\bigl\{j\geq1, F(j)\geq1- 1/n\bigr\}
\\
&=& \inf\biggl\{j\geq1, \sum_{k>j}p_k
\leq1/n\biggr\}.
\end{eqnarray*}
Omitting the slowly varying functions, when $\vec\nu(1/\cdot)\in
\textsc
{RV}_\alpha$, $0<\alpha<1$, $j^*$ is of order $n^{{\alpha}/{(1-\alpha)}}$
and $p_{j^*}$ is of order $n^{-{1}/{(1-\alpha)}}$. The closer to
$1$ is $\alpha$, the smaller the probability of $j^*$. When $\alpha$
goes to $0$, this probability becomes $1/n$. With geometric
frequencies, $j^*$ is $\frac{\log(n)}{\log (1/1-q )}$ and
$p_{j^*}$ is $\frac{q}{n(1-q)}$.\vspace*{2pt} Hence, around the quantile of order
$1-1/n$, there are symbols which may contribute significantly to the
missing mass' fluctuations.

Another interesting case consists of distributions which are very
light-tailed, in the sense that $\frac{p_{k+1}}{p_k}\to0$ when $k\to
\infty$. An example of these is the Poisson distribution $\mathcal
{P}(\lambda)$, for which $\frac{p_{k+1}}{p_k}=\frac{\lambda}{k}
\mathrel{{\rightarrow}_{k\to+\infty}} 0 $. The next proposition
shows that for such concentrated distributions, the missing mass
essentially concentrates on two points.

\begin{prop}\label{prop:potpourri}
In the infinite urns scheme with probability mass function
$(p_k)_{k\in\mathbb{N}}$, if $p_k>0$ for all $k$ and $\lim_{k\to
\infty} \frac{p_{k+1}}{p_k}= 0$, then there exists a sequence of
integers $(u_n)_{n\in\mathbb{N}}$ such that
\[
\lim_{n\to\infty} \PROB \bigl\{ M_{n,0} \in \bigl\{
\overline{F}(u_n),\overline{F}(u_n+1) \bigr\} \bigr\} = 1
,
\]
where $\overline{F}(k)= \sum_{j>k} p_j $.
\end{prop}

\section{Applications}
\label{sec:applications}

\subsection{Estimating the regular variation index}
\label{sec:estim-alpha}

When working in the regular variation setting, the most basic
estimation task is to estimate the regular variation index $\alpha$. We
already mentioned in Section~\ref{sec:notation} the fact that, when
$\vec{\nu} \in\textsc{rv}_\alpha, \alpha\in(0,1)$, the ratio
$K_{n,1}/K_n$ provides a consistent estimate of $\alpha$. This is
actually only one among a family of estimators of $\alpha$ that one may
construct. The next result shows this, and is a direct consequence of
Proposition~\ref{prop:asymp-cumulated-occup}.

\begin{prop}\label{prop:alpha:estimate}
If $\vec{\nu} \in\textsc{rv}_\alpha, \alpha\in(0,1]$, then for all
$r\geq1$
\[
\frac{r K_{n,r}}{K_{n,\overline{r}}}
\]
is a strongly consistent estimator of $\alpha$.
\end{prop}

Thus, writing $k_n=\max{\{r, K_{n,r} >0\}}$, at time $n$, we can have
up to $k_n$ non-trivial estimators of $\alpha$. One would expect these
estimators to offer various bias-variance trade-offs, and one could
ostensibly select an ``optimal'' $r$ via model selection.

\subsection{Estimating the missing mass}
\label{sec:good-turing}

The Good--Turing estimation problem (Good \cite{Good1953}) is that of
estimating $M_{n,r}$ from the observation $(U_1, U_2,\ldots, U_n)$.
For large scores $r$, designing estimators for $M_{n,r}$ is
straightforward: we assume that the empirical distribution mimics the
sampling distribution, and that the empirical probabilities $\frac{r
K_{n,r}}{n}$ are likely to be good estimators. The question is more
delicate for rare events. In particular, for $r=0$, it may be a bad
idea to assume that there is no missing mass $M_{n,0}=0$, that is to
assign a zero probability to the set of symbols that do not appear in
the sample. Various ``smoothing'' techniques have thus developed, in
order to adjust the maximum likelihood estimator and obtain more
accurate probabilities.

In particular, Good--Turing estimators attempt to estimate $(M_{n,r})_r$
from $(K_{n,r})_r$ for all $r$. They are defined as
\[
G_{n,r}=\frac{(r+1)K_{n,r+1}}{n} .
\]

The rationale for this choice comes from the following observations.
%
\begin{equation}
\label{eq:1} \EXP G_{n,0}= \frac{\EXP
 [ K_{n,1} ]}{n}= \EXP
M_{n-1,0} = \EXP M_{n,0} + \frac
{\EXP
M_{n,1}}{n}
\end{equation}
and
%
\begin{equation}
\label{eq:2} \EXP G_{n,r} = \frac{(r+1)\EXP K_{n,r+1}}{n} = \EXP
M_{n-1,r}.
\end{equation}
In the Poisson setting, there is no bias: $\EXP G_r(t)=(r+1)\frac{\EXP
K_{r+1}(t)}{t}=\EXP M_r(t)$.

Here, we primarily focus on the estimation of the missing masses
$M_{n,0}$ and $M_0(t)$, though most of the methodology extends also to
$r>0$, with the appropriate concentration results. From \eqref{eq:1}
and \eqref{eq:2}, Good--Turing estimators look like slightly biased
estimators of the relevant masses. In particular, the bias $\EXP
G_{n,0}-\EXP M_{n,0}$ is always positive but smaller than $1/n$. It is
however far from obvious to determine scenarios where these estimators
are consistent and where meaningful confidence regions can be constructed.

When trying to estimate the missing mass $M_{n,0}$ or $\EXP M_{n,0}$,
consistency needs to be redefined since the estimand is not a fixed
parameter of interest but a random quantity whose expectation further
depends on $n$. Additive consistency, that is bounds on $\widehat
{M}_{n,0}-M_{n,0}$ is not a satisfactory notion, because, as $M_{n,0}$
tends to $0$, the trivial constant estimator $0$ would be universally
asymptotically consistent. Relative consistency, that is control on
$(\widehat{M}_{n,0}-M_{n,0})/M_{n,0}$ looks like a much more reasonable
notion. It is however much harder to establish.

In order to establish relative consistency of a missing mass estimator,
we have to check that $\EXP[\widehat{M}_{n,0}-M_{n,0}]$ is not too
large with respect to $\EXP M_{n,0}$, and that both $\widehat{M}_{n,0}$
and $M_{n,0}$ are concentrated around their mean values.

As shown in Ohannessian and Dahleh \cite{ohannessian2012rare}, the
Good--Turing estimator of
the missing mass is not universally consistent in this sense. This
occurs principally in very light tails, such as those described in
Section~\ref{sec:too-slow}.
%
\begin{prop}[(Ohannessian and Dahleh \cite{ohannessian2012rare})]
When the sampling distribution is geometric with small enough $q\in
(0,1)$, there exists $\eta>0$, and a subsequence $n_i$ such that for
$i$ large enough, $G_{n_i,0}/M_{n_i,0}=0$ with probability no less than
$\eta$.
\end{prop}

On the other hand, the concentration result of Corollary~\ref
{prop:bernstein:right:regvar} gives a law of large numbers for
$M_{n,0}$ (by a direct application of the Borel--Cantelli lemma), which
in turn implies the strong multiplicative consistency of the
Good--Turing estimate.
%
\begin{cor} \label{cor:LGNmissingmass}
We have the following two regimes of consistency for the Good--Turing
estimator of the missing mass.
\begin{longlist}[(ii)]
\item[(i)] If the counting function $\vec{\nu}$ is such that $\EXP
K_{n,2}/\EXP K_{n,1}$ remains bounded and $\EXP K_{n,1}\to+\infty$ (in
particular, when $\vec\nu$ is regularly varying with index $\alpha
\in
(0,1]$ or $\alpha=0$ and $\vec\nu\in\Pi_a$ with $a\to\infty$),
\[
\frac{M_{n,0}}{\EXP M_{n,0}} \overset{\mathbb{P}} {\rightarrow} 1 ,
\]
and the Good--Turing estimator of $M_{n,0}$ defined by
$G_{n,0}=K_{n,1}/n$, is multiplicatively consistent in probability:
\[
\frac{G_{n,0}}{M_{n,0}}\overset{\mathbb{P}} {\rightarrow} 1.
\]

\item[(ii)] If furthermore $\EXP K_{n,\overline{2}} /\EXP K_{n,1}$ remains bounded
and if, for all $c>0$,
$\sum_{n=0}^{\infty} \exp(-c\EXP\times   K_{n,1}) <\infty$ (in
particular, when $\vec\nu$ is regularly varying with index $\alpha
\in
(0,1]$), then these two convergences occur almost surely.
\end{longlist}
\end{cor}

\begin{rem}
One needs to make assumptions on the sampling distribution to guarantee
the consistency of the Good--Turing estimator. In fact, there is no hope
to find a universally consistent estimator of the missing mass without
any such restrictions, as shown recently by Mossel and Ohannessian
\cite
{MosselOhannessian2015}. 
\end{rem}

Consistency is a desirable property, but the concentration inequalities
provide us with more power, in particular in terms of giving confidence
intervals that are asymptotically tight. For brevity, we focus here on
the Poisson setting to derive concentration inequalities which in turn
yield confidence intervals. A similar, but somewhat more tedious,
methodology yields confidence intervals in the binomial setting as well.

\subsubsection{Concentration inequalities for
$G_0(t)-M_0(t)$}
\label{sec:conc-ineq-gtt}

In the Poisson setting, the analysis of the Good--Turing estimator is
illuminating. As noted earlier, the first pleasant observation is that
the Good--Turing estimator is an unbiased estimator of the missing mass.
Second, the variance of $G_0(t) -M_0(t)$ is simply related to
occupancy counts:
%
\begin{equation}
\label{eq:3} \var\bigl(G_0(t) -M_0(t)\bigr) =
\frac{1}{t^2} \bigl(\EXP K_1(t)+ 2 \EXP K_2(t)
\bigr).
\end{equation}
Third, simple yet often tight concentration inequalities can be
obtained for $G_0(t) -M_0(t)$.
%
\begin{prop}
\label{prop:conc-ineq-gtt-1}
The random variable $G_0(t) -M_0(t)$ is sub-gamma on the right tail
with variance factor $\var(G_0(t) -M_0(t))$ and scale factor $1/t$,
and sub-gamma on the left tail with variance factor $3 \EXP K(t)/t^2$
and scale factor $1/t$.

For all $\lambda\geq0$,
\begin{longlist}[(i)]
\item[(i)]
$
\log\EXP\mathe^{\lambda(G_0(t) -M_0(t))} \leq\var(G_0(t)
-M_0(t)) t^2 \phi (\frac{\lambda}{t} ), \textrm{and}
$
\item[(ii)]
$
\log\EXP\mathe^{\lambda(M_0(t)-G_0(t))} \leq\frac{3 \EXP
K(t)}{2 t^2} \frac{\lambda^2}{1-\lambda/t}.
$
\end{longlist}
\end{prop}

We are now in a position to build confidence intervals for the missing mass.
%
\begin{prop}
\label{prop:conf-missing}
With probability larger than $1-4\delta$, the following hold
\begin{eqnarray*}
M_0(t) &\leq& G_0(t) + \frac{1}{t} \biggl(
\sqrt{6 K(t) \log\frac{1}{\delta}} + 5 \log\frac{1}{\delta} \biggr),
\\
M_0(t)& \geq& G_0(t) -\frac{1}{t} \biggl( \sqrt{2
\bigl(K_1(t)+2K_2(t)\bigr) \log\frac{1}{\delta}} + 4
\log\frac{1}{\delta} \biggr).
\end{eqnarray*}
\end{prop}

To see that these confidence bounds are asymptotically tight, consider
the following central limit theorem. A similar results can be
paralleled in the binomial setting.

\begin{prop}\label{prop:TCL-ratio-GT-Mn0}
If the counting function $\vec{\nu}$ is regularly varying with index
$\alpha\in(0,1]$, the following central limit theorem holds for the
ratio $G_0(t)/M_0(t)$:
\[
\frac{\EXP K_1(t)}{\sqrt{\EXP K_1(t)+2\EXP K_2(t)}} \biggl(\frac
{G_0(t)}{M_0(t)} - 1 \biggr) \rightsquigarrow
\mathcal{N}(0,1).
\]
\end{prop}

\begin{rem}
Note that when $\alpha=1$, this convergence occurs faster: the speed is
of order $\sqrt{n\ell_1(n)}$ instead of $\sqrt{n^\alpha\ell(n)}$.
\end{rem}

\subsection{Estimating the number of species}
\label{sec:estim-numb-disc}

Fisher's number of species problem (Fisher \textit{et al.} \cite{Fisher1943})
consists of
estimating $K_{(1+\tau) n}-K_n$ for $\tau>0$, the number of distinct
new species one would observe if the data collection runs for an
additional fraction $\tau$ of time. This was posed primarily within the
Poisson model in the original paper (Fisher \textit{et al.} \cite{Fisher1943})
and later by
Efron and Thisted \cite{efron1976estimating}, but the same question
may also be asked in
the binomial model. The following estimates come from straightforward
computations on the asymptotics given in Theorems \ref{th:portmanteau},
\ref{thmm:portmanteau-rapid-var} and \ref{thmm:portmanteau-slow-var}.

\begin{prop}
If the counting function $\vec{\nu}$ is regularly varying with index
$\alpha\in(0,1]$, letting $\hat{\alpha}$ be any of the estimates
$rK_{n,r}/K_{n,\overline{r}}$ of $\alpha$ from Proposition~\ref
{prop:alpha:estimate}, then any of the following quantities
\[
\bigl(\tau^{\hat{\alpha}}-1\bigr)K_n ,\qquad \frac{\tau^{\hat{\alpha}}-1}{\hat
{\alpha}}K_{n,1}
 \quad\mathrm{and}\quad \Biggl(\prod_{k=2}^r
\frac
{k}{k-1-\hat{\alpha}} \Biggr) \frac{\tau^{\hat{\alpha}}-1}{\hat{\alpha}}K_{n,r} ,\qquad  r\geq2 ,
\]
is a strongly consistent estimate of $K_{\tau n}-K_n$, the number of
newly discovered species when the sample size is multiplied by $\tau$.

If the counting function $\vec{\nu}$ is in $\Pi_{a}$, with $a(n)\to
+\infty$, then, for each $r\geq1$,
\[
\log(\tau)rK_{n,r}
\]
is an estimate of $K_{\tau n}-K_n$, consistent in probability.
\end{prop}

\section{Discussion}
\label{sec:discussion}

To conclude the paper, we review our results in a larger context, and
propose some connections, extensions, and open problems.

\subsection{The cost of Poissonization and negative correlation}
\label{sec:cost-poisson}
Resorting to Poissonization or negative correlation may have a price.
It may lead to variance overestimates.
Gnedin \textit{et al.} (\cite{gnedin2007notes}, Lemma~1), asserts that
for some
constant $c$
\[
\bigl\llvert \var\bigl(K(n)\bigr) -\var(K_n) \bigr\rrvert \leq
\frac{c}{n} \max \bigl( 1 , \EXP K_1(n)^2 \bigr).
\]
This bound conveys a mixed message. As $\EXP K_1(n)/n$ tends to $0$, it
asserts that
\[
\bigl\llvert \var\bigl(K(n)\bigr) -\var(K_n) \bigr\rrvert /\EXP
K_1(n)
\]
tends to $0$. But there exist scenarios where $\EXP K_1(n)^2/n$ tends
to infinity. It is shown in Gnedin \textit{et al.} \cite
{gnedin2007notes} that
$\EXP
K_1(n)^2/(n\var(K(n)))$ tends to $0$, so that, as soon as $n\var(K(n))$
tends to infinity (which might not always be the case), the two
variances $\var(K_n)$ and $\var(K(n))$ are asymptotically equivalent.

It would be interesting to find necessary and sufficient conditions
under which there is equivalence. Though these aren't generally known,
it is instructive to compare $ \var(K(n))$, $\var(K_n)$ and
$\operatorname{Var}^{\mathrm{ind}}
(K_n)$ the variance upper bound obtained from negative correlation by
bounding their differences. For instance, one can show that for any
sampling distribution we have:
\[
\frac{\EXP K_2(2n)}{n} \leq\var\bigl(K(n)\bigr)- \operatorname{Var}^{\mathrm
{ind}}(K_n)
\leq\frac
{2\EXP
K_2(n)}{n}
\]
and
\[
0 \leq\operatorname{Var}^{\mathrm{ind}}(K_n)-\var(K_n)
\leq\frac
{(\EXP K_{n,1})^2}{n} - \frac
{\EXP K_{2n,2}}{2n-1}.
\]
These bounds are insightful but, without any further assumptions on the
sampling distribution, they are not sufficient to prove asymptotic equivalence.

\subsection{Extensions of regular variation}
\label{sec:extensions}

The regular variation hypothesis is an elegant framework, allowing one
to derive, thanks to Karamata and Tauberian theorems, simple and
intelligible equivalents for various moments. As we have seen in Remark~\ref{rem:necessary}, regular variation comes very close to being a
necessary condition for exponential concentration. It may however seem
too stringent. Without getting too specific, let us mention that other
less demanding hypotheses also yield the asymptotic relative orders
that work in favor of the concentration of the missing mass. For
instance, referring back to Remark~\ref{rem:necessary}, one could
instead ask for:
\[
0< \liminf_{t\to\infty} \frac{\Phi_{\overline{2}}(t)}{\Phi
_2(t)}\leq \limsup
_{t\to\infty} \frac{\Phi_{\overline{2}}(t)}{\Phi_2(t)} < \infty .
\]
Recalling that $\Phi_{\overline{2}}'(t)=\frac{2\Phi_2(t)}{t}$, and
applying Corollary~2.6.2. of Bingham \textit{et al.} \cite{BiGoTe87}, one
obtains that $\Phi_2$
is in the class $OR$ of $O$-regularly varying functions and $\Phi
_{\overline{2}}$ is in the class $ER$ of extended regularly varying
functions, that is, for all $\lambda\geq1$
\[
0< \liminf_{t\to\infty} \frac{\Phi_2(\lambda t)}{\Phi_2(t)}\leq \limsup
_{t\to\infty} \frac{\Phi_2(\lambda t)}{\Phi_2(t)} <\infty
\]
and
\[
\lambda^d\leq\liminf_{t\to\infty} \frac{\Phi_{\overline
{2}}(\lambda
t)}{\Phi_{\overline{2}}(t)}\leq
\limsup_{t\to\infty} \frac{\Phi
_{\overline{2}}(\lambda t)}{\Phi_{\overline{2}}(t)} \leq\lambda^c ,
\]
for some constants $c$ and $d$. Observe that this result, which is the
equivalent of Karamata's theorem, differs from the regular variation
setting, in the sense that the control on the derivative $\Phi_2$ is
looser than the one on $\Phi_{\overline{2}}$, whereas, in the Karamata
theorem, both the function and its derivative inherit the regular
variation property.

We can in turn show that $\Phi(t)=\EXP K(t)$ is in the class $OR$ and,
by Theorem~2.10.2 of Bingham \textit{et al.} \cite{BiGoTe87}, this is
equivalent to $\vec\nu
(1/\cdot) \in OR$, as $\Phi$ is the Laplace--Stieltjes transform of
$\vec
\nu$.

\subsection{Random measures}

As noted by Gnedin \textit{et al.} \cite{gnedin2007notes}, the
asymptotics for
the moments of
the occupancy counts in the regular variation setting is still valid
when the frequencies $(p_j)_{j\geq1}$ are random, in which case the
measure $\nu$ is defined by
\[
\EXP \Biggl[\sum_{j=1}^\infty
f(p_j) \Biggr]= \int_0^1 f(x)\nu
(\mathrm {d}x) ,
\]
for all functions $f\geq0$. We can also define the measure $\nu_1$ by
\[
\EXP \Biggl[\sum_{j=1}^\infty
p_jf(p_j) \Biggr]= \int_0^1
f(x)\nu _1(\mathrm{d}x) ,
\]
for all functions $f\geq0$. This measure corresponds to the
distribution of the frequency of the first discovered symbol.

For instance, when $(p_j)_{j\geq1}$ are Poisson--Dirichlet($\alpha
$,$0$) with $0<\alpha<1$, the measure $\nu_1(\mathrm{d}x)$ is the
size-biased distribution of $\operatorname{PD}(\alpha,0)$, that is Beta $\mathcal
{B}(1-\alpha,\alpha)$ (see Pitman and Yor \cite{pitman1997two}).
Thus, we have:
\begin{eqnarray*}
\nu_1[0,x]&=&\frac{1}{\mathcal{B}(1-\alpha,\alpha)} \int_0^x
t^{-\alpha
}(1-t)^{\alpha-1}\,\mathrm{d}t
\\
&\mathop{\sim}\limits_{x\to0}& \frac{x^{1-\alpha}}{(1-\alpha)\mathcal
{B}(1-\alpha,\alpha)}
\end{eqnarray*}
and, by Gnedin \textit{et al.} \cite
{gnedin2007notes}, Proposition~13, this
is equivalent to
\[
\vec{\nu}(x)\mathop{\sim}\limits_{x\to0}\frac{1}{\alpha\mathcal{B}(1-\alpha
,\alpha)}x^{-\alpha}.
\]

Thus, denoting by $N(x)$ the random number of frequencies $p_j$ which
are larger than $x$, the expectation $\vec\nu(x)=\EXP N(x)$ is
regularly varying. One can also show that the mass-partition mechanism
of the distribution $\operatorname{PD}(\alpha, 0)$ almost surely generates $N(x)$ to
be regularly varying. To see this, refer to Pitman and Yor \cite{pitman1997two},
Proposition~10
or to Bertoin
\cite{bertoin2006random}, Proposition~2.6,
which assert that the limit
\[
L:= \lim_{n\to\infty} np_n^\alpha
\]
exists almost surely. This is equivalent to
\[
N(x)\mathop{\sim}\limits_{x\to0} x^{-\alpha}L\qquad\mbox{almost surely.}
\]

The $\operatorname{PD}(\alpha,0)$ distribution can be generated through a Poisson
process with intensity measure $\nu([x,\infty])=cx^{-\alpha}$. Without
entering into further details, let us mention that similar almost sure
results hold even when the intensity measure $\nu$ is not a strict
power, but satisfies the property
\[
\nu\bigl([x,\infty]\bigr)\mathop{\sim}\limits_{x\to0} x^{-\alpha}\ell(x) ,
\]
with $\ell$ slowly varying, Gnedin \cite{MR2684164}, Section~6. Working
with a regular variation hypothesis thus gives us more flexibility than
assuming specific Bayesian priors.

\section{Proofs}
\label{sec:proofs}

\subsection{Fundamental techniques}

\subsubsection{Efron--Stein--Steele inequalities}
\label{sec:efron-stein}

Our variance bounds mostly follow from the Efron--Stein--Steele
inequality (Efron and Stein~\cite{efron:stein:1981}), which states that
when a random
variable is expressed as a function of many independent random
variables, its variance can be controlled by the sum of the local fluctuations.

\begin{thmm}\label{the:es:ineq}
Let $\mathcal{X}$ be some set, $(X_1,X_2,\ldots,X_n)$ be independent
random variables taking values in $\mathcal{X}$, $f : \mathcal{X}^n
\rightarrow\R$ be a measurable function of $n$ variables, and
$Z=f(X_1,X_2,\ldots,X_n)$.

For all $i\in\{1,\ldots,n\}$, let $X^{(i)}=(X_1,\ldots
,X_{i-1},X_{i+1},\ldots,X_n)$ and $\EXP^{(i)}Z=\EXP[Z|X^{(i)}]$. Then,
letting $v=\sum_{i=1}^n \EXP[(Z-\EXP^{(i)}Z)^2]$,
\[
\var[Z]\leq v.
\]

If $X_{1}',\ldots,X_{n}'$ are independent copies of $X_1,\ldots,X_n$,
then letting $Z_i'=f(X_1,\ldots,X_{i-1},\break X_i', X_{i+1},\ldots,X_n)$,
\[
v= \sum_{i=1}^n \EXP\bigl[
\bigl(Z-Z_i'\bigr)_{+}^2\bigr]
\leq\sum_{i=1}^n \EXP\bigl[(Z-Z_i)^2
\bigr] ,
\]
where the random variables $Z_i$ are arbitrary $X^{(i)}$-measurable and
square-integrable random variables.
\end{thmm}
%

\subsubsection{Negative association}
\label{sec:negative-association}

The random variables $K_n$, $K_{n,r}$, and $M_{n,r}$ are sums or
weighted sums of Bernoulli random variables. These summands depend on
the scores $(X_{n,j})_{j\geq1}$ and therefore are not independent.
Transforming the fixed-$n$ binomial setting into a continuous time
Poisson setting is one way to circumvent this problem. This is the
Poissonization method. In this setting, the score variables
$(X_j(n))_{j\geq1}$ are independent Poisson variables with respective
means $np_j$. Results valid for the Poisson setting can then be
transferred to the fixed-$n$ setting, up to approximation costs. For
instance, Gnedin \textit{et al.} \cite{gnedin2007notes} (Lemma~1) provide
bounds on the
discrepancy between expectations and variances in the two settings.
(See also our discussion in Section~\ref{sec:cost-poisson}.)

Another approach to deal with the dependence is to invoke the notion of
negative association, which provides a systematic comparison between
moments of certain monotonic functions of the occupancy scores. In our
present setting, this will primarily be useful for bounding the
logarithmic moment generating function, which is an expectation of
products, by products of expectations, thus recovering the structure of
independence. This has already been used to derive exponential
concentration for occupancy counts
(see Dubhashi and Ranjan \cite{dubhashi:ranjan:1998},
Shao \cite{shao2000comparison},
McAllester and Schapire~\cite{McAllesterSchapire2000},
Ohannessian and Dahleh \cite{ohannessian2012rare}).
It is also useful for bounding variances. We use this notion throughout
the proofs, and therefore present it here formally.

\begin{dfn}[(Negative association)]
Real-valued random variables $Z_1,\ldots, Z_K$ are said to be negatively
associated if, for any two disjoint subsets $A$ and $B$ of $\{1,\ldots
,K\}$, and any two real-valued functions $f:\mathbb{R}^{|A|}\mapsto
\mathbb{R}$ and $g:\mathbb{R}^{|B|}\mapsto\mathbb{R}$ that are both
either coordinate-wise non-increasing or coordinate-wise
non-decreasing, we have:
\[
\EXP \bigl[f(Z_A) \cdot g(Z_B) \bigr] \leq\EXP
\bigl[f(Z_A) \bigr]\cdot \EXP \bigl[g(Z_B) \bigr].
\]
\end{dfn}

In particular, as far as concentration properties are concerned, sums
of negatively associated variables can only do better than sums of
independent variables.

\begin{thmm}[(Dubhashi and Ranjan \cite{dubhashi:ranjan:1998})]
For each $n\in\mathbb{N}$, the occupancy scores $(X_{n,j})_{j\geq1}$
are negatively associated.
\end{thmm}

As monotonic functions of negatively associated variables are also
negatively associated, the variables $(\IND_{\{X_{n,j}>0\}})_{j\geq1}$
(respectively, $(\IND_{\{X_{n,j}=0\}})_{j\geq1}$) are negatively
associated as increasing (respectively, decreasing) functions of
$(X_{n,j})_{j \geq1}$. This is of pivotal importance for our proofs of
concentration results for $K_n$ and $M_{n,0}$. For $r\geq1$, the
variables $(\IND_{\{X_{n,j}= r\}})_{j\geq1}$ appearing in $K_{n,r}$
are not negatively associated. However, following Ohannessian and
Dahleh \cite
{ohannessian2012rare}, one way to deal with this problem is to observe
that
\[
K_{n,r}= K_{n,\overline{r}}-K_{n,\overline{r+1}} ,
\]
recalling that $K_{n,\overline{r}}=\sum_{j=1}^\infty\IND_{\{
X_{n,j}\geq r\}}$ is the number of urns that contain at least $r$ balls
and that the Bernoulli variables appearing in $K_{n,\overline{r}}$ are
negatively associated.

\subsubsection{Potter's inequalities and its variants}
\label{sec:potter}

One useful result from regular variation theory is provided by Potter's
inequality (see Bingham \textit{et al.} \cite{BiGoTe87}, {de Haan}
and Ferreira \cite{HaaFei06}, for proofs and
refinements).

\begin{thmm}[(Potter--Drees inequalities)]
\label{thmm:potter:rv}
\begin{longlist}[(ii)]
\item[(i)] If $f \in\textsc{rv}_\gamma$, then for all $\delta>0$, there
exists $t_0=t_0(\delta)$,
such that for all $t,x\colon\min(t,tx)>t_0$,
\[
(1-\delta) x^\gamma\min \bigl(x^\delta,x^{-\delta} \bigr)
\leq \frac
{f(tx)}{f(t)} \leq(1+\delta) x^\gamma\max \bigl(x^\delta,x^{-\delta
}
\bigr).
\]

\item[(ii)] If $\ell\in\Pi_a$, then for all $\delta_1, \delta_2$, there
exists $t_0$ such that for all $t\geq t_0$, for all $x\geq1$,
\[
(1-\delta_2)\frac{1-x^{\delta_1}}{\delta_1}-\delta_2<
\frac{\ell
(tx)-\ell(t)}{a(t)}< (1+\delta_2)\frac{x^{\delta_1}-1}{\delta
_1}+\delta
_2.
\]
\end{longlist}
\end{thmm}
%

\subsection{Occupancy counts}

\subsubsection{Variance bounds for occupancy counts}
\label{sec:proof-vari-bounds-occ-counts}

%
\begin{pf*}{Proof of Proposition~\ref{prop:varKt}}
Recall that in the Poisson setting,
\[
\frac{\mathrm{d}\EXP K(t)}{\mathrm{d}t} = \frac{\mathbb
{E}K_1(t)}{t} = {\mathbb{E}M_0(t)}.
\]
This entails
\begin{eqnarray*}
\var\bigl(K(t)\bigr)= \sum_{j=1}^\infty
\mathe^{-tp_j}\bigl(1-\mathe^{-tp_j}\bigr)
= \EXP K(2t)-\EXP K(t)
= \int_t^{2t} \EXP M_0(s)
\,\mathrm{d}s.
\end{eqnarray*}
Now, as $\EXP M_0(s) $ is non-increasing,
\[
\frac{\EXP K_1(2t)}{2}=t \EXP M_0(2t) \leq\var\bigl(K(t)\bigr)\leq t
\EXP M_0(t) = \EXP K_1(t).
\]

Moving on to the binomial setting, let $K_n^i$ denote the number of
occupied urns when the $i$th ball is replaced by an
independent copy. Then
\[
\var(K_n) \leq\EXP \Biggl[\sum_{i=1}^n
\bigl(K_n -K_n^{i}\bigr)_+^2
\Biggr] ,
\]
where $(K_n -K_n^{i})_+$ denotes the positive part. Now,
$K_n-K_{n}^{i}$ is positive if and only if ball $i$ is moved from a
singleton into in a nonempty urn. Thus $\var(K_n)\leq\EXP
[K_{n,1}(1-M_{n,0}) ]$.
\end{pf*}

%
\begin{pf*}{Proof of Proposition~\ref{prop:varKnrinf}}
The bound $r\EXP K_{n,r}$ follows from the Efron--Stein inequality:
denoting by $K_{n,\overline{r}}^{(i)}$ the number of cells with occupancy
score larger than $r$ when ball $i$ is removed, then
\[
K_{n,\overline{r}}-K_{n,\overline{r}}^{(i)} = %
\cases{ 1, &\quad $
\mbox{if ball $i$ is in a $r$-ton},$ \vspace*{2pt}
\cr
0, &\quad $\mbox{otherwise.}$}
\]
And thus, we get $\sum_{i=1}^n
(K_{n,\overline{r}}-K_{n,\overline{r}}^{(i)})^2=rK_{n,r}$.

The second bound follows from the negative association of the variables
$(\IND_{\{X_{n,j}\geq r\}})_j$ (negative correlation is actually sufficient):
\[
\var \Biggl(\sum_{j=1}^\infty
\IND_{\{X_{n,j}\geq r\}} \Biggr) \leq \sum_{j=1}^\infty
\var(\IND_{\{X_{n,j}\geq r\}}) \leq\EXP K_{n,\overline{r}} .
\]
\upqed\end{pf*}

\subsubsection{Concentration inequalities for occupancy counts}
\label{sec:proof-conc-ineq-occ-counts}

%
\begin{pf*}{Proof of Proposition~\ref{prop:conc-ineq-knrinf}}
Let $X_{n,j}$ denote the occupancy score of cell $j, j\in\mathbb{N}$, then
\[
K_{n,\overline{r}} = \sum_{j=1}^\infty
\mathbb{I}_{\{X_{n,j}\geq r\}}.
\]
As noted in Section~\ref{sec:negative-association}, $K_{n,\overline
{r}}$ is
a sum of negatively associated Bernoulli random variables. Moreover,
the Efron--Stein inequality implies that for each $j\in\mathbb{N}$,
\[
\var(\mathbb{I}_{\{X_{n,j}\geq r\}}) \leq r \EXP\mathbb{I}_{\{
X_{n,j}= r\}}.
\]

Thus we have
\begin{eqnarray*}
\log\EXP\mathe^{\lambda(K_{n,\overline{r}}-\EXP K_{n,\overline
{r}})} &\leq &\sum_{j=1}^\infty
\log\EXP\mathe^{\lambda(\mathbb{I}_{\{
X_{n,j}\geq r\}} -\EXP\mathbb{I}_{\{X_{n,j}\geq r\}})}
\\
& \leq& \sum_{j=1}^\infty\var(
\mathbb{I}_{\{X_{n,j}\geq r\}}) \phi (\lambda)
\\
& \leq& \phi(\lambda) \sum_{j=1}^\infty r
\EXP\IND_{\{X_{n,j}=r\}
}
\\
& = & \phi(\lambda) r \EXP K_{n,r} ,
\end{eqnarray*}
where the first inequality comes from negative association, the second
inequality is Bennett's inequality for Bernoulli random variables, and
the last inequality comes from the Efron--Stein inequality.
The other bound comes from the fact that $\var(\mathbb{I}_{\{
X_{n,j}\geq r\}}) \leq\EXP\mathbb{I}_{\{X_{n,j}\geq r\}}$.
\end{pf*}

%
\begin{pf*}{Proof of Proposition~\ref{prop:conc-ineq-knr}}
As $K_{n,r}=K_{n,\overline{r}}-K_{n,\overline{r+1}}$,
\begin{eqnarray*}
&&\{ K_{n,r} \geq\EXP K_{n,r} + x \}
\\
&&\quad\subseteq \biggl\{ K_{n,\overline{r}} \geq\EXP K_{n,\overline{r}} +
\frac
{x}{2} \biggr\} \cup \biggl\{K_{n,\overline{r+1}} \leq\EXP
K_{n,\overline{r+1}} -\frac{x}{2} \biggr\}.
\end{eqnarray*}
By Proposition~\ref{prop:conc-ineq-knrinf}, Bernstein
inequalities hold for both $K_{n,\overline{r}}$ and $K_{n,\overline
{r+1}}$, with
variance proxies $\EXP rK_{n,r}$ (or $\EXP K_{n,\overline{r}}$) and
$(r+1)K_{n,r+1}$ (or $\EXP K_{n,\overline{r+1}}\leq\EXP
K_{n,\overline{r}}$)
respectively. Hence,
\begin{eqnarray*}
&&\mathbb{P} \{ K_{n,r} \geq\EXP K_{n,r} +{x} \}
\\
&&\quad \leq \exp \biggl( - \frac{x^2/4}{2(r\EXP K_{n,r} + x/6)} \biggr) + \exp \biggl( -
\frac{x^2/4}{2((r+1)\EXP K_{n,r+1} )} \biggr)
\\
& &\quad \leq 2 \exp \biggl( - \frac{x^2/4}{2(\max(r\EXP K_{n,r},(r+1)
K_{n,r+1}) + x/6)} \biggr).
\end{eqnarray*}
The same reasoning works for the alternative variance proxies and for
the left tails.
\end{pf*}

\subsection{Missing mass}
\label{sec:proof-missing-mass}

\subsubsection{Variance bounds for the missing mass}
\label{sec:proof-variance-missing-mass}

%
\begin{pf*}{Proof of Proposition~\ref{prop:vari-missingmass}}
In the Poisson setting,
\begin{eqnarray*}
\var\bigl(M_0(t)\bigr)= \sum_{j=1}^\infty
p_j^2 \mathe^{-tp_j} \bigl(1-
\mathe^{-tp_j} \bigr) \leq\sum_{j=1}^\infty
p_j^2 \mathe^{-tp_j} = \frac{2}{t^2} \EXP
K_2(t).
\end{eqnarray*}
In the binomial setting, by negative correlation,
\[
\var(M_{n,0}) \leq\sum_{j=1}^\infty
p_j^2 \bigl(1- (1-p_j)^n \bigr)
(1-p_j)^n \leq\sum_{j=1}^\infty
p_j^2 \mathe^{-n p_j} = \frac{2}{n^2} \EXP
K_2(n).
\]
\upqed\end{pf*}

\subsubsection{Concentration inequalities for the missing mass}
\label{sec:proof-conc-ineq-missing-mass}

%
\begin{pf*}{Proof of Proposition~\ref{prop:bennett:left:mn0}}
For all $\lambda\in\mathbb{R}$,
\begin{eqnarray*}
\log\EXP \bigl[ \mathe^{\lambda(M_{n,0}-\EXP M_{n,0})} \bigr] &=& \log\EXP \bigl[
\mathe^{\lambda\sum_{j=1}^\infty p_j
(Y_j-\EXP Y_j)} \bigr]
\\
&\leq& \sum_{j=1}^\infty\log\EXP \bigl[
\mathe^{\lambda p_j (Y_j -
\EXP[Y_j])} \bigr]
\\
&\leq& \sum_{j=1}^\infty(1-p_j)^n
\bigl(1-(1-p_j)^n\bigr) \phi(\lambda p_j) ,
\end{eqnarray*}
where the first inequality comes from negative association, and the
second is Bennett's inequality for Bernoulli random variables.

Noting that $\lim_{\lambda\to0_-} \phi(\lambda)/\lambda^2 =
\lim_{\lambda\to0_+}\phi(\lambda)/\lambda^2 = 1/2$, the function
$\lambda\mapsto
{\phi(\lambda)}/{\lambda^2}$ has a continuous increasing
extension on $\mathbb{R}$.
Hence, 
for $\lambda\leq0$, we have $\phi(\lambda)\leq{\lambda^2}/{2}$.

Thus, for $\lambda<0$,
\begin{eqnarray*}
\log\EXP \bigl[ \mathe^{\lambda(M_{n,0}-\EXP M_{n,0})} \bigr] &\leq& \sum
_{j=1}^\infty p_j^2
(1-p_j)^n\bigl(1-(1-p_j)^n
\bigr) \frac{\lambda^2}{2}
\\
&=& \sum_{j=1}^\infty p_j^2
\var[Y_j] \frac{\lambda^2}{2}.
\end{eqnarray*}
Recall that $\sum_{j=1}^\infty p_j^2 \var[Y_j] \leq2 \EXP K_2(n) /n^2$
(Proposition~\ref{prop:vari-missingmass}).
\end{pf*}

%
\begin{pf*}{Proof of Proposition~\ref{prop:log-laplace-missing-mass}}
From the beginning of the proof of Proposition~\ref
{prop:bennett:left:mn0}, that is, thanks to negative association and to
the fact that each Bernoulli random variable satisfies a Bennett inequality,
\[
\log\EXP \bigl[ \mathe^{\lambda(M_{n,0}-\EXP M_{n,0})} \bigr] 
\leq \sum
_{j=1}^\infty\mathe^{-np_j} \phi(
\lambda p_j).
\]
Now, using the power series expansion of $\phi$,
\begin{eqnarray*}
\sum_{j=1}^\infty\mathe^{-np_j} \phi(
\lambda p_j)&=& \sum_{j=1}^\infty
\mathe^{-np_j} \sum_{r=2}^\infty
\frac{(\lambda p_j)^r}{r!}
\\
&=& \sum_{r=2}^\infty \biggl(
\frac{\lambda}{n} \biggr)^r \sum_{j=1}^\infty
\mathe^{-np_j} \frac{(np_j)^r}{r!}.
\end{eqnarray*}
We recognize that for each $r$, $\sum_{j=1}^\infty
\mathe^{-np_j} \frac{(np_j)^r}{r!} =\EXP K_r(n)$, so that
\begin{eqnarray*}
\log\EXP \bigl[ \mathe^{\lambda(M_{n,0}-\EXP M_{n,0})} \bigr]&\leq& \sum
_{r=2}^\infty \biggl(\frac{\lambda}{n}
\biggr)^r\EXP K_r(n).
\end{eqnarray*}
\upqed\end{pf*}

%
\begin{pf*}{Proof of Theorem~\ref{thmm:struct-ineq-right}}
Using Proposition~\ref{prop:log-laplace-missing-mass} and noticing
that for each $r\geq2$,
$\EXP K_r(n)\leq\EXP K_{\overline{2}}(n)$, we immediately obtain that
\begin{eqnarray*}
\log\EXP \bigl[ \mathe^{\lambda(M_{n,0}-\EXP M_{n,0})} \bigr]&\leq& \EXP K_{\overline{2}}(n)
\sum_{r=2}^\infty \biggl(\frac{\lambda
}{n}
\biggr)^r
\\
& = & \lambda^2 \frac{ \EXP K_{\overline{2}}(n)/n^2}{1-\lambda/n} ,
\end{eqnarray*}
which concludes the proof.
\end{pf*}

\subsection{Regular variation}

%
\begin{pf*}{Proof of Proposition~\ref{prop:asymp-cumulated-occup}}
By monotonicity of $K_{n,\overline{r}}$, we have the following strong
law for any sampling distribution
\[
K_{n,\overline{r}}=\sum_{s=r}^\infty
K_{n,s}\mathop{\sim}\limits_{+\infty}\sum_{s=r}^\infty
\EXP K_s(n)\qquad \mbox{a.s.},
\]
(see Gnedin \textit{et al.} \cite{gnedin2007notes}, the discussion
after Proposition~2).
Recall that $X_j(n)\sim\mathcal{P}(np_j)$ and that, if $Y\sim
\mathcal
{P}(\lambda)$, then $\PROB [Y\leq k ] = \frac{\Gamma
(k+1,\lambda)}{k!}$, where $\Gamma(z,x)=\int_x^{+\infty} \mathe
^{-t}t^{z-1}\,\mathrm{d}t$ is the incomplete Gamma function. Hence
\begin{eqnarray*}
\sum_{s=r}^\infty\EXP K_s(n)&=&
\sum_{j= 1}^\infty\PROB \bigl[X_j(n)
\geq r \bigr]
\\
&=& \sum_{j= 1}^\infty\frac{1}{(r-1)!}\int
_0^{np_j} \mathe^{-t}t^{r-1}
\,\mathrm{d}t
\\
&=& \frac{1}{(r-1)!} \int_0^1 \int
_0^{nx} \mathe^{-t}t^{r-1}
\mathrm {dt}\cdot\nu(\mathrm{d}x)
\\
&=& \frac{1}{(r-1)!} \int_0^1 n
\mathe^{-nx}(nx)^{r-1}\vec\nu(x)\, \mathrm {d}x
\\
&=& \frac{1}{(r-1)!}\int_0^{+\infty}
\mathe^{-z}z^{r-1}\vec\nu (z/n)\,\mathrm{d}z
\\
&\mathop{\sim}\limits_{+\infty}& \frac{\vec\nu(1/n)}{(r-1)!}\Gamma(r-\alpha).
\end{eqnarray*}

In particular,
\[
K_{n,\overline{r}} \mathop{\sim}\limits_{+\infty}\frac{rK_{n,r}}{\alpha} \qquad\mbox
{a.s.}
\]
\upqed\end{pf*}

%
\begin{pf*}{Proof of Theorem~\ref{thmm:bernstein:right:slowvar}}
Let us recall Proposition~\ref{prop:log-laplace-missing-mass}:
\[
\log\EXP \bigl[ \mathrm{e}^{\lambda(M_{n,0}-\EXP M_{n,0})} \bigr] \leq\sum
_{r=2}^\infty \biggl(\frac{\lambda}{n}
\biggr)^r \EXP K_r(n).
\]

Now, bounding each $\EXP K_r(n)$ by $\EXP K_{\overline{2}}(n)$ is not
sufficient to get the right order for the variance: $\EXP K_{\overline
{2}}(n)$ is of order $\ell(n)$ whereas $\var M_{n,0}$ is of order $a(n)/n^2$.

We explore more carefully the structure of $\EXP K_r(n)$ and show that
these quantities are uniformly (in $r$) bounded by a function of order
$a(n)$ for large enough $n$, that is, that there exists $n^*\in\N$ and
$C\in\R_+$ such that for all $n\geq n^*$, for all $r\geq1$, $\EXP
K_r(n)\leq Ca(n)$.

Before going into the proof, we observe that for $r\geq n/a(n)$, the
result is true. Indeed, from the identity $\sum_{r=1}^\infty r\EXP
K_r(n)=n$, we deduce that $r\EXP K_{\overline{r}}(n)\leq n$, so that
for $r\geq n/a(n)$, $\EXP K_{\overline{r}}(n)\leq a(n)$. Thus, we
assume that $r\leq n/a(n)$.

First, we easily deal with the contribution to $\EXP K_r(n)$ of the
symbols with probability less than $1/n$. Indeed
\begin{eqnarray*}
I_1^r&:=&\int_0^{1/n}
\mathe^{-nx}\frac{(nx)^r}{r!}\nu(\mathrm{d}x) \leq\int
_0^{1/n} \mathe^{-nx}\frac{(nx)^2}{2!}
\nu(\mathrm{d}x) \leq\EXP K_2(n).
\end{eqnarray*}

As $\EXP K_2(n)\sim a(n)/2$, for all $\delta_0$, there exists $n_0$
such that for all $n\geq n_0$, for all $r\geq1$, $I_{1}^r \leq
(1+\delta_0)/2$.

For the contribution of the symbols with probability larger than
$1/n$, integration by part and change of variable yield:
\begin{eqnarray*}
I_2^r&:=& \int_{1/n}^1
\mathe^{-nx}\frac{(nx)^r}{r!}\nu(\mathrm {d}x)
\\
&=& \biggl[\mathe^{-nx}\frac{(nx)^r}{r!}\bigl(-\vec\nu(x)\bigr)
\biggr]_{1/n}^1 + \int_{1/n}^1
\mathe^{-nx} \frac{n^r}{r!}\bigl(rx^{r-1}-nx^r
\bigr)\vec\nu (x)\,\mathrm {d}x
\\
&=& \frac{\vec\nu(1/n)\mathrm{e}^{-1}}{r!} + \int_1^\infty
\mathe^{-z} \biggl(\frac
{z^{r-1}}{(r-1)!}-\frac{z^r}{r!} \biggr)\vec
\nu(z/n)\,\mathrm{d}z.
\end{eqnarray*}

As $\int_1^\infty\mathe^{-z} (\frac{z^{r-1}}{(r-1)!}-\frac
{z^r}{r!} )\,\mathrm{d}z=-\PROB[\mathcal{P}(1)=r]=-\mathrm{e}^{-1}/r!$, we
can rearrange the previous expression:
\[
I_2^r=a(n)\int_1^\infty
\mathe^{-z} \biggl(\frac
{z^{r-1}}{(r-1)!}-\frac
{z^r}{r!} \biggr)
\frac{\vec\nu(z/n)-\vec\nu(1/n)}{a(n)}\,\mathrm {d}z.
\]

Notice that when $z\in[1,r]$, the integrand is negative, so we simply
ignore this part of the integral and restrict ourselves to
\[
I_3^r:=\int_r^\infty
\mathe^{-z} \biggl(\frac{z^{r}}{r!}-\frac
{z^{r-1}}{(r-1)!} \biggr)
\frac{\vec\nu(1/n)-\vec\nu
(z/n)}{a(n)}\,\mathrm {d}z,
\]
which we try to bound by a constant term for $n$ greater than some
integer that does not depend on $r$.

The main ingredient of our proof is the next version of the
Potter--Drees inequality (see Theorem~\ref{thmm:potter:rv} in
Section~\ref{sec:potter}
and {de Haan} and Ferreira \cite{HaaFei06}, point $4$ of Corollary
B.2.15): for $\ell\in\Pi_a$, for arbitrary $\delta_1,
\delta_2$, 
there exists $t_0$ such that for all $t\geq t_0$, and for all $x\leq1$
with $tx\geq t_0$,
\[
(1-\delta_2)\frac{1-x^{-\delta_1}}{\delta_1}-\delta_2<
\frac{\ell
(t)-\ell(tx)}{a(t)}< (1+\delta_2)\frac{x^{-\delta_1}-1}{\delta
_1}+\delta
_2.
\]

Thus, for arbitrary $\delta_1$, $\delta_2$, there exists $n_1$ such
that, for all $n\geq n_1$, for all $z \in[1,n/n_1]$,
\[
\frac{\vec\nu(1/n)-\vec\nu(z/n)}{a(n)} \leq(1+\delta_2)\frac
{z^{\delta
_1}-1}{\delta_1}+
\delta_2.
\]

As $r\leq n/a(n)$, taking, if necessary, $n$ large enough so that
$a(n)\geq n_1$, we have $r\leq n/n_1$ and
\begin{eqnarray*}
I_{3}^r &\leq& \int_r^{n/n_1}
\mathe^{-z} \biggl(\frac
{z^r}{r!}-\frac
{z^{r-1}}{(r-1)!} \biggr)
\biggl((1+\delta_2)\frac{z^{\delta
_1}-1}{\delta
_1}+\delta_2 \biggr)
\,\mathrm{d}z
\\
&&{} + \int_{n/n_1}^\infty\mathe^{-z} \biggl(
\frac
{z^r}{r!}-\frac{z^{r-1}}{(r-1)!} \biggr)\frac{\vec\nu(1/n)-\vec\nu
(z/n)}{a(n)}\,\mathrm{d}z
\\
&=:& I_{4}^r + I_{5}^r ,
\end{eqnarray*}
with
\begin{eqnarray*}
I_{4}^r&\leq& \delta_2 + \frac{1+\delta_2}{\delta_1}
\int_r^\infty\mathe^{-z} \biggl(
\frac{z^{r+\delta_1}}{r!}-\frac
{z^r}{r!}+\frac{z^{r-1}}{(r-1)!}-\frac{z^{r-1+\delta
_1}}{(r-1)!}
\biggr)\,\mathrm{d}z
\\
&\leq& \delta_2 + \frac{1+\delta_2}{\delta_1} \int_r^\infty
\mathe^{-z} \biggl(\frac{z^{r+\delta_1}}{r!}- \frac
{z^{r-1+\delta_1}}{(r-1)!} \biggr)
\,\mathrm{d}z
\\
&=& \delta_2 + \frac{1+\delta_2}{\delta_1} \biggl(\frac{\Gamma
(r+1+\delta
_1 ,r)}{\Gamma(r+1)}-
\frac{\Gamma(r+\delta_1 ,r)}{\Gamma(r)} \biggr) ,
\end{eqnarray*}
where $\Gamma(a,x)=\int_x^\infty\mathe^{-t}t^{a-1} \,\mathrm{d}t$ is the
incomplete Gamma function. Using the fact that $\Gamma
(a,x)=(a-1)\Gamma
(a-1,x)+x^{a-1}\mathe^{-x}$, we have
\begin{eqnarray*}
I_{4}^r &\leq& \delta_2 +
\frac{1+\delta_2}{\delta_1\Gamma(r+1)} \bigl(\Gamma(r+1+\delta_1 ,r)- (r+
\delta_1)\Gamma(r+\delta_1 ,r)+ \delta_1
\Gamma(r+\delta_1 ,r) \bigr)
\\
&=& \delta_2 + \frac{1+\delta_2}{\delta_1} \biggl(\frac{r^{r+\delta
_1}\mathe^{-r}}{r!}+
\delta_1\frac{\Gamma(r+\delta_1 ,r)}{\Gamma(r+1)} \biggr).
\end{eqnarray*}

By Stirling's inequality, for all $r$,
\[
\frac{r^{r+\delta_1}\mathe^{-r}}{r!}\leq r^{\delta_1}(2\pi r)^{-1/2}.
\]
Thus, taking $\delta_1 = 1/4$, the right-hand term is uniformly bounded
by $1$. And $\frac{\Gamma(r+\delta_1 ,r)}{\Gamma(r+1)}$ is also bounded
by 1. Thus
\[
I_{4}^r \leq\delta_2 +
\frac{1+\delta_2}{\delta_1}(1+\delta_1)
\]
and
\begin{eqnarray*}
I_{5}^r&\leq& \frac{\vec\nu(1/n)}{a(n)}\int
_{n/n_1}^\infty\mathe ^{-z} \biggl(
\frac{z^r}{r!}-\frac{z^{r-1}}{(r-1)!} \biggr) \,\mathrm {d}z
\\
&=& \frac{\vec\nu(1/n)}{a(n)} \mathe^{-n/n_1} \frac{(n/n_1)^r}{r!}
\\
&\leq&\frac{\vec\nu(1/n)}{a(n)} \mathe^{-n/n_1} \frac
{(n/n_1)^{\lfloor
n/n_1 \rfloor}}{\lfloor n/n_1 \rfloor!}.
\end{eqnarray*}
By Stirling's inequality, this bound is smaller than $\frac{\vec\nu
(1/n)}{a(n)}(2\pi(n/n_1))^{-1/2}$,
which tends to $0$ as $n\to\infty$. Thus, there exists $n_2$ such that
for all $n\geq n_2$, and all $r\leq n/n_1$, $I_{5}^r\leq\delta_2$.

In the end, we get that for all $\delta_0\geq0$, $\delta_1$ with
$0\leq\delta_1\leq1/4$, and $\delta_2\geq0$, there exists
$n^*=\max
(n_0,n_1,n_2)$ such that for all $n\geq n^*$, for all $r\geq1$,
\[
\EXP K_r(n) \leq a(n) \biggl( \frac{1+\delta_0}{2} +
\delta_2 +\frac
{1+\delta_2}{\delta_1}(1+\delta_1) +
\delta_2 \biggr).
\]

Taking for instance $\delta_1=1/4$ and $\delta_0=\delta_2=1/15$, we
have that for large enough $n$ and for all $r\geq1$,
\[
\EXP K_r(n) \leq6a(n)
\]
and
\[
\log\EXP \bigl[ \mathrm{e}^{\lambda(M_{n,0}-\EXP M_{n,0})} \bigr] \leq\frac{12a(n)}{n^2}\cdot
\frac{\lambda^2}{2(1-\lambda/n)}.
\]
\upqed\end{pf*}

%
\begin{pf*}{Proof of Proposition~\ref{prop:potpourri}}
Under the condition of the Proposition~\ref{prop:potpourri}, from Gr{\"
u}bel and Hitczenko \cite
{MR2582705}, with probability tending to $1$, the sample is gap-free,
hence the missing mass is $\overline{F}(\max(X_1,\ldots,X_n))$.

The condition of the proposition implies the condition described in
Anderson \cite{And70}, i.e. $\lim_{n\to+\infty} \frac{\overline
{F}(n+1)}{\overline{F}(n)}=0$, to ensure the existence of a sequence of
integers $(u_n)_{n\in\mathbb{N}}$ such that
\[
\lim_{n\to\infty} \PROB \bigl\{ \max(X_1,
\ldots,X_n) \in \{u_n,u_n+1 \} \bigr\} = 1.
\]
\upqed\end{pf*}

\subsection{Applications}

%
\begin{pf*}{Proof of Corollary~\ref{cor:LGNmissingmass}}
Let us assume that $\EXP K_{n,1}\to\infty$. Using the fact that
$0\leq
\EXP G_{n,0} - \EXP M_{n,0}\leq1/n$, we notice that as soon as $\EXP
K_{n,1}\to\infty$, $\EXP G_{n,0}\sim_{+\infty}\EXP M_{n,0}$.
Now by Chebyshev's inequality,
\begin{eqnarray*}
\PROB \biggl[ \biggl\llvert \frac{M_{n,0}}{\EXP M_{n,0}}-1 \biggr\rrvert >\varepsilon
\biggr] &\leq&\frac{\var(M_{n,0})}{\varepsilon^2 (\EXP
M_{n,0})^2} \leq \frac{2\EXP K_2(n)}{\varepsilon^2n^2(\EXP M_{n,0})^2}
\\
&\sim& \frac{2(\EXP K_{n,2} +o(1))}{\varepsilon^2 (\EXP K_{n,1})^2} ,
\end{eqnarray*}
where we used that $\vert\EXP K_2(n) - \EXP K_{n,2}\vert\to0$ (see
Lemma~1, Gnedin
\textit{et al.} \cite{gnedin2007notes}). On the other hand,
\[
\PROB \biggl[ \biggl\llvert \frac{K_{n,1}}{\EXP K_{n,1}}-1 \biggr\rrvert >\varepsilon
\biggr] \leq \frac{\var(K_{n,1})}{\varepsilon^2
(\EXP K_{n,1})^2} \leq\frac{\EXP K_{n,1}+2\EXP K_{n,2}}{\varepsilon^2(\EXP K_{n,1})^2} ,
\]
showing that if, furthermore, $\EXP K_{n,2}/\EXP K_{n,1}$ remains
bounded, the ratios $M_{n,0}/\EXP M_{n,0}$, $G_{n,0}/ \EXP G_{n,0}$ and
thus $M_{n,0}/G_{n,0}$ converge to $1$ in probability.
To get almost sure convergence, we use Theorem~\ref
{thmm:struct-ineq-right} to get that when $\EXP K_{n,1}\to\infty$,
\begin{eqnarray*}
\PROB \biggl[ \biggl\llvert \frac{M_{n,0}}{\EXP M_{n,0}}-1 \biggr\rrvert >\varepsilon
\biggr] &\leq& 2\exp \biggl(-\frac{\varepsilon^2(\EXP
M_{n,0})^2}{2(2\EXP K_{\overline{2}}(n)/n^2 + \EXP M_{n,0}/n)} \biggr)
\\
&= & 2\exp \biggl(-\frac{\varepsilon^2(\EXP K_{n,1}+o(\EXP
K_{n,1}))^2}{2(2\EXP K_{n,\overline{2}} + \EXP K_{n,1} +o(\EXP
K_{n,1}))} \biggr).
\end{eqnarray*}
If $\EXP K_{n,\overline{2}}/\EXP K_{n,1}$ remains bounded, this becomes
smaller than $c_1\exp (-c_2\varepsilon^2 \EXP K_{n,1} )$. Hence,
if $\exp(-c\EXP K_{n,1})$ is summable for all $c>0$, we can apply the
Borel--Cantelli lemma and obtain the almost sure convergence of
$M_{n,0}/\EXP M_{n,0}$ to $1$. Moreover, by Proposition~\ref
{prop:conc-ineq-knr},
\begin{eqnarray*}
\PROB \biggl[ \biggl\llvert \frac{K_{n,1}}{\EXP K_{n,1}}-1 \biggr\rrvert >\varepsilon
\biggr] &\leq& 4\exp \biggl(-\frac{\varepsilon^2(\EXP
K_{n,1})^2}{2(4\max(\EXP K_{n,1}, 2\EXP K_{n,2})+2/3)} \biggr) ,
\end{eqnarray*}
which shows that under these assumptions $K_{n,1}/\EXP K_{n,1}$ also
tends to $1$ almost surely.
\end{pf*}

%
\begin{pf*}{Proof of Proposition~\ref{prop:conc-ineq-gtt-1}}
The random variable $G_0(t)-M_0(t)$ is a sum of independent, centered
and bounded random variables, namely
\[
G_0(t)-M_0(t) = \frac{1}{t} \sum
_{j=1}^\infty\IND_{X_j(t)=1} - tp_j
\IND_{X_j(t)=0}.
\]

Bound (i) follows immediately from the observation that each $\IND_{X_j(t)=1}
- tp_j \IND_{X_j(t)=0}$ satisfies a
Bennett inequality,
\begin{eqnarray*}
\log\EXP\mathe^{\lambda(G_0(t) -M_0(t))} &\leq& \sum_{j=1}^\infty
\var(\IND_{X_j(t)=1} - tp_j \IND_{X_j(t)=0}) \phi \biggl(
\frac{\lambda}{t} \biggr)
\\
& = & \var\bigl(G_0(t) -M_0(t)\bigr) t^2
\phi \biggl(\frac{\lambda}{t} \biggr).
\end{eqnarray*}

Bound (ii) follows from the observation that each $ \IND
_{X_j(t)=0}-\frac{1}{tp_j}\IND_{X_j(t)=1}$ satisfies a Bennett inequality,
\begin{eqnarray*}
\log\EXP\mathe^{\lambda(M_0(t)-G_0(t))} &\leq& \sum_{j=1}^\infty
\var \biggl(\IND_{X_j(t)=0}-\frac{1}{tp_j}\IND_{X_j(t)=1} \biggr)
\phi ( {\lambda p_j} )
\\
& = & \sum_{j=1}^\infty \biggl(1+
\frac{1}{tp_j} \biggr) \mathe^{-tp_j} \phi(\lambda p_j)
\\
& = & \sum_{r\geq2} \biggl(\frac{\lambda}{t}
\biggr)^r \sum_{j=1}^\infty
\biggl(1+\frac{1}{tp_j} \biggr)\mathe^{-tp_j} \frac{(tp_j)^r}{r!}
\\
& = & \sum_{r\geq2} \biggl(\frac{\lambda}{t}
\biggr)^r \biggl( \EXP K_r(t) +\frac{1}{r}\EXP
K_{r-1}(t) \biggr)
\\
& \leq& \sum_{r\geq2} \biggl(\frac{\lambda}{t}
\biggr)^r \frac{3
\EXP
K(t)}{2} ,
\end{eqnarray*}
which concludes the proof.
\end{pf*}

%
\begin{pf*}{Proof of Proposition~\ref{prop:conf-missing}}
With probability greater than $1-2\delta$, by Proposition~\ref
{prop:conc-ineq-gtt-1},
\[
G_0(t) -M_0(t) \leq\frac{1}{t}\sqrt{2 \bigl(
\EXP K_1(t) +2 \EXP K_2(t)\bigr)\log \frac{1}{\delta}} +
\frac{\log({1}/{\delta})}{3t}
\]
and
\[
G_0(t) -M_0(t) \geq- \frac{1}{t}\sqrt{6 \EXP
K(t)\log\frac{1}{\delta}} - \frac{\log({1}/{\delta})}{t}.
\]
We may now invoke concentration inequalities for $K_1(t)+2 K_2(t) $ and
$K(t)$. Indeed, with probability greater than $1-\delta$,
$K(t)\geq\EXP K(t) -\sqrt{2 \EXP K(t) \log\frac{1}{\delta}}$ which
entails $\sqrt{\EXP K(t)}\leq\sqrt{K(t)+\frac{\log
({1}/{\delta})}{2}}+\sqrt{\frac{\log({1}/{\delta})}{2}}$.

We have $2K_2(t)+K_1(t) \geq2\EXP K_2(t)+\EXP
K_1(t)-\sqrt{ 4 (2\EXP K_2(t)+\EXP
K_1(t)) \log\frac{1}{\delta}} $ with probability greater than
$1-\delta
$, which entails
\[
\sqrt{2\EXP K_2(t)+\EXP K_1(t)}\leq\sqrt{
\bigl(2K_2(t)+K_1(t) \bigr)+ \log\frac{1}{\delta}} +
\sqrt{{\log\frac{1}{\delta}}} ,
\]
which concludes the proof.
\end{pf*}

%
\begin{pf*}{Proof of Proposition~\ref{prop:TCL-ratio-GT-Mn0}}
The covariance matrix $\cov(t)$ of $(G_0(t),M_0(t))$ can be written in
terms of the expected occupancy counts as
\[
\cov(t) = \frac{1}{t^2} %
\pmatrix{ \EXP K_1(t) & 0
\vspace*{2pt}
\cr
0 & 2\EXP K_2(t) } %
- \frac{\EXP K_2(2t)}{2t^2}
\pmatrix{ 1 & 1\vspace*{2pt}
\cr
1 & 1 } %
.
\]

From Karlin \cite{Kar67}, we have
\[
\cov(t)^{-1/2} %
\pmatrix{ G_0(t) - \EXP
G_0(t) \vspace*{2pt}
\cr
M_0(t) - \EXP
M_0(t) } %
\rightsquigarrow\mathcal{N}(0,I_2),
\]
where $I_2$ is the identity matrix, which can be rewritten as
\[
\Sigma(t)^{-1/2} %
\pmatrix{ \displaystyle\frac{G_0(t)}{\EXP G_0(t)}-1
\vspace*{2pt}
\cr
\displaystyle\frac{M_O(t)}{\EXP M_0(t)}-1} %
\rightsquigarrow
\mathcal{N}(0,I_2) ,
\]
with $\Sigma(t)= (\EXP G_0(t))^{-2} \cov(t)$.

The delta method applied to the function $(x_1,x_2)\mapsto x_1/x_2$ yields
\[
\biggl((1\ -1) \Sigma(t) %
\pmatrix{ 1 \vspace*{2pt}
\cr
-1 } %
 \biggr)^{-1/2} \biggl(\frac{G_0(t)}{M_0(t)} - 1 \biggr) \rightsquigarrow
\mathcal{N}(0,1)
\]
and
\[
\biggl((1\ -1) \Sigma(t) %
\pmatrix{ 1 \vspace*{2pt}
\cr
-1 } %
 \biggr)^{-1/2}= \frac{\EXP K_1(t)}{\sqrt{\EXP K_1(t)+2\EXP K_2(t)}} ,
\]
which concludes the proof.
\end{pf*}

\begin{rem}
The proof for the binomial setting is very similar, the only difficulty
being that $\EXP G_{n,0}$ and $\EXP M_{n,0}$ are no longer equal.
However, the bias becomes negligible with respect to the fluctuations,
that is, for $v_n$ either $n^\alpha\ell(n)$ or $n\ell_1(n)$
\[
\sqrt{v_n} \biggl(\frac{\EXP
G_{n,0}}{\EXP M_{n,0}}-1 \biggr)\mathop{\to}_{n\to\infty}0.
\ \ 
\]
\end{rem}


%
%

\printhistory
\end{document}